\documentclass{amsart}

\calclayout

\usepackage{amssymb,amscd}
\usepackage{verbatim}
\usepackage{lamsarrow}

\usepackage{array,dcolumn}
\newcolumntype{C}{>{$}c<{$}}
\newcolumntype{L}{>{$}l<{$}}
\newcolumntype{R}{>{$}r<{$}}

\def\to{\mathchoice
{\longrightarrow}
{\rightarrow}
{\rightarrow}
{\rightarrow}}

\def\mapsto{\DOTSB\mapstochar\to}

\def\hookrightarrow{\mathchoice
{\DOTSB\lhook\joinrel\relbar\joinrel\rightarrow}
{\DOTSB\lhook\joinrel\rightarrow}
{\DOTSB\lhook\joinrel\rightarrow}
{\DOTSB\lhook\joinrel\rightarrow}}

\newtheorem{theorem}{Theorem}
\newtheorem*{theorem*}{Theorem}
\newtheorem{proposition}[theorem]{Proposition}
\newtheorem{lemma}[theorem]{Lemma}

\theoremstyle{definition}

\newtheorem{definition}[theorem]{Definition}

\newcommand{\Z}{\mathbb{Z}}
\newcommand{\C}{\mathbb{C}}
\newcommand{\R}{\mathbb{R}}
\newcommand{\Q}{\mathbb{Q}}

\newcommand{\om}{\omega}
\renewcommand{\o}{\otimes}
\newcommand{\p}{\partial}
\newcommand{\bull}{\bullet}

\newcommand{\<}{\langle}
\renewcommand{\>}{\rangle}
\renewcommand{\]}{{]\!]}}
\renewcommand{\[}{{[\![}}

\newcommand{\CM}{\mathcal{M}}
\newcommand{\Mbar}{\overline{\mathcal{M}}}
\newcommand{\Cbar}{\overline{\mathcal{C}}}
\renewcommand{\SS}{\mathbb{S}}
\newcommand{\VV}{\mathbb{V}}
\renewcommand{\P}{\mathbb{P}}
\DeclareMathOperator{\Serre}{\mathsf{e}}
\DeclareMathOperator{\SP}{Sp}
\DeclareMathOperator{\SL}{SL}

\newcommand{\EE}{\mathbb{E}}
\DeclareMathOperator{\Aut}{Aut}
\DeclareMathOperator{\Tr}{Tr}

\newcommand{\Om}{\Omega}
\newcommand{\MM}{\mathsf{M}}
\newcommand{\tMM}{\widetilde{\mathsf{M}}}
\newcommand{\MMbar}{\overline{\mathsf{M}}}
\newcommand{\CCbar}{\overline{\mathsf{C}}}
\newcommand{\virt}{{\textup{virt}}}

\newcommand{\stab}{\mathsf{stab}}

\newcommand{\Nov}{\Lambda}
\DeclareMathOperator{\NN}{N}
\DeclareMathOperator{\NE}{NE}
\DeclareMathOperator{\ZZ}{Z}
\DeclareMathOperator{\ZE}{ZE}
\DeclareMathOperator{\supp}{supp}
\DeclareMathOperator{\ev}{ev}
\DeclareMathOperator{\Edge}{Edge}
\DeclareMathOperator{\Flag}{Flag}
\DeclareMathOperator{\VERT}{Vert}

\DeclareMathOperator{\Res}{Res}
\newcommand{\Dual}{\vee}
\newcommand{\CO}{\mathcal{O}}
\newcommand{\GG}{\mathbb{G}}
\newcommand{\CA}{\mathcal{A}}
\newcommand{\Abar}{\overline{\mathcal{A}}}
\newcommand{\CJ}{\mathcal{J}}
\newcommand{\Jbar}{\overline{\mathcal{J}}}
\DeclareMathOperator{\diag}{diag}
\DeclareMathOperator{\GSP}{GSp}

\DeclareMathOperator{\GL}{GL}
\newcommand{\CW}{\mathcal{W}}
\newcommand{\CE}{\mathcal{E}}
\DeclareMathOperator{\cd}{cd}

\begin{document}

\title{Topological Recursion Relations in Genus $2$}

\author{E. Getzler}

\address{Department of Mathematics, Northwestern University, Evanston,
IL 60208, USA}

\email{getzler@math.nwu.edu}

\maketitle

Consider a two-dimensional topological field theory coupled to topological
gravity with primary fields $\{\gamma_a\}$ (see Witten \cite{Witten} for an
explanation of this notation). An example of such a theory is provided by
the Gromov-Witten invariants of a smooth projective variety $V$; in this
case, the primary fields are a basis of the cohomology $H^\bull(V,\Q)$.

The metric on the space of primaries is defined by the formula
$$
\eta_{ab} = \<P\gamma_a\gamma_b\>_0 ,
$$
where $P=\gamma_0$ is the puncture operator. (In the case of Gromov-Witten
invariants, $\gamma_0$ is the identity of $H^\bull(V,\Q)$.) Denote by
$\eta^{ab}$ the metric on the dual space, and by
$\gamma^a=\eta^{ab}\gamma_b$ the dual basis.

Let $\{t^a_k\mid k\ge0\}$ be the coordinates on the large phase
space. Introduce generating functions for the correlation functions, by the
formula
\begin{multline} \label{<<>>}
\<\<\tau_{k_1}(x_1)\dots\tau_{k_n}(x_n)\>\>_g = \Bigl\langle
\tau_{k_1}(x_1)\dots\tau_{k_n}(x_n) \exp \sum_{k,a} t_k^a \tau_k(\gamma_a)
\Bigr\rangle_g \\ \equiv \sum_{N=0}^\infty \frac{1}{N!}
\sum_{\substack{l_1\dots l_N \\ a_1 \dots a_N}} t_{l_N}^{a_N} \dots
t_{l_1}^{a_1} \< \tau_{k_1}(x_1)\dots\tau_{k_n}(x_n)
\tau_{l_1}(\gamma_{a_1}) \dots \tau_{l_N}(\gamma_{a_N}) \>_g .
\end{multline}

We can now formulate the topological recursion relations satisfied by
the genus $0$ correlation functions:
\begin{equation} \label{trr0}
\<\<\tau_{k_1}(x_1)\tau_{k_2}(x_2)\tau_{k_3}(x_3)\>\>_0 = \sum_a
\<\<\tau_{k_1-1}(x_1)\gamma^a\>\>_0
\<\<\gamma_a\tau_{k_2}(x_2)\tau_{k_3}(x_3)\>\>_0 .
\end{equation}
In combination with other properties of these models, these recursion
relations allow the genus $0$ correlation functions for gravitational
descendents to be expressed in terms of correlation functions for primary
fields. Ultimately, \eqref{trr0} is a consequence of the vanishing of the
tautological classes $\psi_i$ on the Deligne-Mumford compactification
$\Mbar_{0,3}$, which follows from the zero-di\-men\-sio\-nal-ity of
$\Mbar_{0,3}$.

Likewise, the genus $1$ correlation functions satisfy the topological
recursion relations
\begin{equation} \label{trr1}
\<\<\tau_k(x)\>\>_1 = \sum_a \<\<\tau_{k-1}(x)\gamma_a\>\>_0
\<\<\gamma^a\>\>_1 + \frac{1}{24} \sum_a
\<\<\tau_{k-1}(x)\gamma_a\gamma^a\>\>_0 .
\end{equation}
These recursion relations follows from the equality of the tautological
class $\psi_1$ on $\Mbar_{1,1}$ with a boundary class of degree $1/24$ on
$\Mbar_{1,1}$.

The first rigourous construction of gravitational descendents was obtained
by Ruan and Tian \cite{RT}, for the Gromov-Witten invariants of
semi-positive symplectic manifolds. In Part 1 of this paper, we study
gravitational descendents for general projective manifolds, by applying the
Behrend-Fantechi construction \cite{B,BF} of the virtual fundamental
classes on $\Mbar_{g,n}(V,\beta)$ (based on the ideas of Li and Tian
\cite{LT}). We verify the above topological recursion relations in this
setting, along with all of the other formulas contained in Dijkgraaf and
Witten \cite{Witten, DW}. It is useful to read our paper in conjunction
with that of Kontsevich and Manin \cite{KM}; although both discuss similar
questions, the emphasis is quite different.

In Part 2, we calculate the topological recursion relations in genus
$2$. To this end, we first calculate the Hodge polynomial of $\Mbar_{2,2}$:
$$
\sum_{i,j} u^i v^j \dim H^{i,j}(\Mbar_{2,2}) =
1+6uv+14u^2v^2+14u^3v^3+6u^4v^4+u^5v^5 ;
$$
the proof uses the results on mixed Hodge modules proved in \cite{II} and a
number of other miscellaneous results. We then use intersection theory in
$\Mbar_{2,2}$ to find topological recursion relations for genus $2$
correlation functions. The situation is very different from that in lower
genus, since the restrictions of the tautological classes $\psi_i$ to
$\CM_{2,n}$ do not vanish. On the other hand, we will see that quadratic
polynomials in the classes $\psi_i$ are boundary classes on $\Mbar_{2,n}$.%
\footnote{We conjecture that this pattern continues: polynomials of degree
$g$ in the tautological classes $\psi_i$ are boundary classes on
$\Mbar_{g,n}$.}

\raggedbottom

Recall Mumford's definition of the cohomology classes $\lambda_i$ on
$\Mbar_g$. Let $\pi:\Mbar_{g,1}\to\Mbar_g$ be the universal curve over
$\Mbar_g$, and consider the Chern classes $\lambda_i=c_i(\EE)$ of the Hodge
bundle $\EE=\pi_*\om_{\Mbar_{g,1}/\Mbar_g}$. We denote the pullbacks of
$\lambda_i$ to $\Mbar_{g,n}$ by the same symbol. By the formulas
$$
\sum_{i=0}^{2k} (-1)^i \lambda_i\lambda_{2k-i} = 0 ,
$$
the classes $\lambda_{2k}$ are universal polynomials in the classes
$\lambda_{2i+1}$; for example, $\lambda_2=\lambda_1^2/2$. In genus $2$, the
divisor $\lambda_1$ is a boundary class (this is not true in higher
genera), and Mumford \cite{Mumford} proves the formula
$$
\lambda_1 = \frac{1}{10} \bigl( 2 \delta_1 + \delta_0 \bigr) .
$$
He also obtains an explicit expression for $\psi_1^2$ as a boundary cycle
on $\Mbar_{2,1}$, which may be written, following Faber \cite{Faber}, as
\begin{align} \notag
\psi_1^2 & = \begin{picture}(40,25)(60,775)
\put( 80,780){\line( 1, 1){ 13}}
\put( 80,780){\line(-1, 1){ 13}}
\put( 80,780){\line( 0,-1){ 18}}
\put( 66,796){\circle{5}}
\put( 94,796){\circle{5}}
\end{picture}
+ \psi_1\*\lambda_1 - \lambda_2 \\ \label{psi^2}
&= \frac{7}{5}
\begin{picture}(40,35)(60,775)
\put( 80,780){\line( 1, 1){ 13}}
\put( 80,780){\line(-1, 1){ 13}}
\put( 80,780){\line( 0,-1){ 18}}
\put( 66,796){\circle{5}}
\put( 94,796){\circle{5}}
\end{picture}
+ \frac{1}{5}
\begin{picture}(30,35)(65,757)
\put( 80,782){\circle{5}}
\put( 80,770){\circle{20}}
\put( 80,760){\line(0,-1){18}}
\end{picture}
- \frac{1}{120}
\begin{picture}(30,35)(65,757)
\put( 80,776){\circle{20}}
\put( 80,766){\line( 0,-1){9}}
\put( 80,755){\circle{5}}
\put( 80,753){\line( 0,-1){11}}
\end{picture}
+ \frac{13}{120}
\begin{picture}(30,35)(65,757)
\put( 80,770){\circle{20}}
\put( 80,760){\line( 3,-4){10}}
\put( 80,760){\line(-3,-4){12}}
\put( 91,745){\circle{5}}
\end{picture}
+ \frac{1}{120}
\begin{picture}(50,35)(55,760)
\put( 70,765){\circle{20}}
\put( 90,765){\circle{20}}
\put( 80,765){\line(0,-1){20}}
\end{picture} \\[20pt] \notag
\end{align}

The main result of Part 2 is the following explicit formula for
$\psi_1\psi_2$:
\begin{gather} \notag
\psi_1\psi_2 =
3 \begin{picture}(30,35)(65,742)
\put( 80,770){\circle{6}}
\put( 78.5,768.5){$\scriptscriptstyle2$}
\put( 80,747){\lamsvector( 0, 1){ 20}}
\put( 80,747){\line(-2,-3){ 10}}
\put( 80,747){\line( 2,-3){ 10}}
\put( 67,722){$1$}
\put( 88,722){$2$}
\end{picture}
+ \frac{13}{5} \begin{picture}(35,35)(62,747)
\put( 80,750){\line(-3,-4){ 10}}
\put( 80,750){\line( 3,-4){ 10}}
\put( 80,750){\line(-3, 4){ 10}}
\put( 80,750){\line( 3, 4){ 10}}
\put( 69,766){\circle{5}}
\put( 91,766){\circle{5}}
\put( 67,727){$1$}
\put( 88,727){$2$}
\end{picture}
+ \frac{4}{5} \begin{picture}(40,35)(60,742)
\put( 77,793){$1$}
\put( 80,772){\line( 0, 1){ 18}}
\put( 80,770){\circle{5}}
\put( 80,767){\line( 0,-1){ 20}}
\put( 80,747){\line(-2,-3){ 10}}
\put( 80,747){\line( 2,-3){ 10}}
\put( 67,722){$2$}
\put( 91,729){\circle{5}}
\end{picture}
+ \frac{4}{5} \begin{picture}(40,35)(60,742)
\put( 77,793){$2$}
\put( 80,772){\line( 0, 1){ 18}}
\put( 80,770){\circle{5}}
\put( 80,767){\line( 0,-1){ 20}}
\put( 80,747){\line(-2,-3){ 10}}
\put( 80,747){\line( 2,-3){ 10}}
\put( 67,722){$1$}
\put( 90,729){\circle{5}}
\end{picture}
- \frac{4}{5} \begin{picture}(40,60)(60,742)
\put( 80,795){\circle{5}}
\put( 80,772){\line( 0, 1){ 20}}
\put( 80,770){\circle{5}}
\put( 80,767){\line( 0,-1){ 20}}
\put( 80,747){\line(-2,-3){ 10}}
\put( 80,747){\line( 2,-3){ 10}}
\put( 67,722){$1$}
\put( 88,722){$2$}
\end{picture} \\[40pt] \label{M22}
+ \frac{23}{120} \begin{picture}(35,35)(62,750)
\put( 80,765){\circle{20}}
\put( 80,755){\line(-3,-4){12}}
\put( 80,755){\line( 0,-1){15}}
\put( 80,755){\line( 3,-4){10}}
\put( 67,730){$1$}
\put( 77,730){$2$}
\put( 91,739){\circle{5}}
\end{picture}
+ \frac{5}{120} \begin{picture}(30,35)(125,740)
\put(135,785){\line( 0,-1){ 37}}
\put(145,768){\circle{20}}
\put(135,745){\circle{5}}
\put(135,742){\line( 0,-1){ 12}}
\put(133,720){$1$}
\put(133,787){$2$}
\end{picture}
+ \frac{5}{120} \begin{picture}(30,35)(125,740)
\put(135,785){\line( 0,-1){ 37}}
\put(145,768){\circle{20}}
\put(135,745){\circle{5}}
\put(135,742){\line( 0,-1){ 12}}
\put(133,720){$2$}
\put(133,787){$1$}
\end{picture}
- \frac{3}{120} \begin{picture}(45,35)(118,740)
\put(135,775){\circle{20}}
\put(135,765){\line( 0,-1){ 17}}
\put(135,745){\circle{5}}
\put(133,743){\line(-3,-4){ 10}}
\put(137,743){\line( 3,-4){ 10}}
\put(120,720){$1$}
\put(145,720){$2$}
\end{picture} \\[15pt] \notag
+ \frac{7}{15} \begin{picture}(35,35)(63,762)
\put( 80,765){\circle{20}}
\put( 80,777){\circle{5}}
\put( 80,755){\line(-3,-4){10}}
\put( 80,755){\line(3,-4){10}}
\put( 67,730){$1$}
\put( 87,730){$2$}
\end{picture}
+ \frac{1}{15} \begin{picture}(35,35)(63,762)
\put( 80,765){\circle{20}}
\put( 80,777){\circle{5}}
\put( 80,755){\line(0,-1){15}}
\put( 80,779){\line(0,1){12}}
\put( 77,730){$2$}
\put( 77,792){$1$}
\end{picture}
+ \frac{1}{15} \begin{picture}(35,35)(63,762)
\put( 80,765){\circle{20}}
\put( 80,777){\circle{5}}
\put( 80,755){\line(0,-1){15}}
\put( 80,779){\line(0,1){12}}
\put( 77,730){$1$}
\put( 77,792){$2$}
\end{picture}
- \frac{1}{15} \begin{picture}(35,35)(118,752)
\put(135,775){\circle{20}}
\put(135,763){\circle{5}}
\put(135,760){\line( 0,-1){ 17}}
\put(135,743){\line(-3,-4){ 10}}
\put(135,743){\line( 3,-4){ 10}}
\put(122,720){$1$}
\put(143,720){$2$}
\end{picture}
- \frac{1}{72} \begin{picture}(50,35)(45,762)
\put( 60,765){\circle{20}}
\put( 80,765){\circle{20}}
\put( 70,740){\line(0,1){50}}
\put( 67,730){$1$}
\put( 67,792){$2$}
\end{picture} \\[15pt] \notag
\end{gather}
In this article, we content ourselves with proving that \eqref{M22} is a
homological equivalence, although we expect that it is actually a rational
equivalence. A key step in its proof is the fact that the rational
cohomology of $\Mbar_{2,2}$ turns out to be generated as a ring by $6$
tautological divisors; we prove this using Faber's computer program
\cite{Faber:Divisors}, which enables us to calculate the intersection
numbers of quintic monomials in the divisors on $\Mbar_{2,2}$.

Finally, applying the results of Part 1, \eqref{psi^2} and \eqref{M22}
imply respectively the topological recursion relations
\begin{multline} \label{trr2}
\<\<\tau_{k+2}(x)\>\>_2 = \<\<\tau_{k+1}(x)\gamma_a\>\>_0
\<\<\gamma^a\>\>_2 + \<\<\tau_k(x)\gamma_a\>\>_0
\<\<\tau_1(\gamma^a)\>\>_2 \\
\begin{aligned}
{}& - \<\<\tau_k(x)\gamma_a\>\>_0 \<\<\gamma^a\gamma_b\>\>_0
\<\<\gamma^b\>\>_2 + \tfrac{7}{10} \<\<\tau_k(x)\gamma_a\gamma_b\>\>_0
\<\<\gamma^a\>\>_1 \<\<\gamma^b\>\>_1 \\ {}& + \tfrac{1}{10}
\<\<\tau_k(x)\gamma_a\gamma_b\>\>_0 \<\<\gamma^a\gamma^b\>\>_1 -
\tfrac{1}{240} \<\<\tau_k(x)\gamma_a\>\>_1
\<\<\gamma^a\gamma_b\gamma^b\>\>_0 \\ {}& + \tfrac{13}{240}
\<\<\tau_k(x)\gamma_a\gamma^a\gamma_b\>\>_0 \<\<\gamma^b\>\>_1 +
\tfrac{1}{960} \<\<\tau_k(x)\gamma_a\gamma^a\gamma_b\gamma^b\>\>_0 ,
\quad \text{for $k\ge0$.} \end{aligned}
\end{multline}
and
\begin{multline} \label{trr22}
\<\<\tau_{k+1}(x)\tau_{l+1}(y)\>\>_2 =
\<\<\tau_{k+1}(x)\gamma^a\>\>_2 \<\<\tau_l(y)\gamma_a\>\>_0 \\
\begin{aligned}
{}& + \<\<\tau_k(x)\gamma_a\>\>_0 \<\<\tau_{l+1}(y)\gamma^a\>\>_2
- \<\<\tau_k(x)\gamma_a\>\>_0 \<\<\tau_l(y)\gamma_b\>\>_0
\<\<\gamma^a\gamma^b\>\>_2 \\ {}& + 3
\<\<\tau_k(x)\tau_l(y)\gamma_a\>\>_0\<\<\gamma^a\>\>_2 +
\tfrac{13}{10} \<\<\tau_k(x)\tau_l(y)\gamma_a\gamma_b\>\>_0
\<\<\gamma^a\>\>_1 \<\<\gamma^b\>\>_1 \\ {}& + \tfrac{2}{5}
\<\<\tau_k(x)\gamma_a\>\>_1 \<\<\tau_l(y)\gamma^a\gamma_b\>\>_0
\<\<\gamma^b\>\>_1 + \tfrac{2}{5}
\<\<\tau_k(x)\gamma^a\gamma_b\>\>_0 \<\<\tau_l(y)\gamma_a\>\>_1
\<\<\gamma^b\>\>_1 \\ {}& - \tfrac{2}{5} \<\<\tau_k(x)\tau_l(y)\gamma_a\>\>_1
\<\<\gamma^a\gamma_b\>\>_0 \<\<\gamma^b\>\>_1 + \tfrac{23}{240}
\<\<\tau_k(x)\tau_l(y)\gamma_a\gamma^a\gamma_b\>\>_0
\<\<\gamma^b\>\>_1 \\ {}& + \tfrac{1}{48}
\<\<\tau_k(x)\gamma_a\gamma^a\gamma_b\>\>_0
\<\<\tau_l(y)\gamma^b\>\>_1 + \tfrac{1}{48}
\<\<\tau_k(x)\gamma^b\>\>_1
\<\<\tau_l(y)\gamma_a\gamma^a\gamma_b\>\>_0 \\ {}& - \tfrac{1}{80}
\<\<\tau_k(x)\tau_l(y)\gamma^b\>\>_1
\<\<\gamma_a\gamma^a\gamma_b\>\>_0 + \tfrac{7}{30}
\<\<\tau_k(x)\tau_l(y)\gamma_a\gamma_b\>\>_0
\<\<\gamma^a\gamma^b\>\>_1 \\ {}& + \tfrac{1}{30}
\<\<\tau_k(x)\gamma_a\gamma_b\>\>_0
\<\<\tau_l(y)\gamma^a\gamma^b\>\>_1 + \tfrac{1}{30}
\<\<\tau_k(x)\gamma^a\gamma^b\>\>_1
\<\<\tau_l(y)\gamma_a\gamma_b\>\>_0 \\ {}& + \tfrac{1}{30}
\<\<\tau_k(x)\tau_l(y)\gamma_a\>\>_0
\<\<\gamma^a\gamma_b\gamma^b\>\>_1 + \tfrac{1}{576}
\<\<\tau_k(x)\tau_l(y)\gamma_a\gamma^a\gamma_b\gamma^b\>\>_0 ,
\quad k,l\ge0.
\end{aligned}
\end{multline}

\subsection*{Acknowledgements}
This paper owes a great debt to C. Faber, both for many extremely helpful
conversations, and for generously sharing the computer program
\cite{Faber:Divisors}, which played an essential role both in the proofs
and in the research behind this paper. I am grateful to M. Kapranov for
telling me about Faltings's Eichler spectral sequence in higher genus, and
also to P. Belorousski, R. Hain and G.~Mess for their helpful remarks.

I thank the Mittag-Leffler Institute, and its director D. Laksov, for
an invitation which provided the opportunity to begin this research,
and the Centre for Mathematics and its Applications of the ANU, and
G. Lehrer, for providing the setting in which it was completed.

The author is partially supported by the NSF.

\newpage

\flushbottom

\part{Gravitational descendents}

\section{Dual graphs}

\subsection{Dual graphs of prestable curves}

\begin{definition}
A pointed prestable curve $(C,z_1,\dots,z_n)$ is a projective, connected
curve $C$ (over $\C$) with no singularities other than simple double
points, together with $n$ distinct, smooth marked points $(z_1,\dots,z_n)$.
\end{definition}

Given a pointed prestable curve $C=(C,z_1,\dots,z_n)$, we will define a
graph $C(G)$, the dual graph of $C$, which represents the homeomorphism
type of $C$. The graph $G=G(C)$ has one vertex $v$ for each component
$C(v)$ of (the normalization of) the curve $C$, labelled by the genus
$g(v)$ of this component.

The edges of the dual graph $G(C)$ correspond to double points of the curve
$C$; the two ends of an edge are attached to the vertices associated to the
components on which the two branches of the double point lie. (If both
branches lie in the same component of $C$, then the edge is a loop.)

Finally, to each marked point $z_i$ of the curve corresponds a leg of the
graph, labelled by $i$, at the vertex corresponding to the component of $C$
on which $z_i$ lies (which is uniquely determined, since $z_i$ is a smooth
point).

In drawing dual graphs, we denote vertices of genus $0$ either by a solid
circle $\begin{picture}(10,8)(0,0)\put(5,3){\circle*{4}}\end{picture}$ or
leave them unmarked, vertices of genus $1$ by an open circle
$\begin{picture}(10,8)(0,0)\put(5,3){\circle{4}}\end{picture}$, and
vertices of genus $g>1$ by $\begin{picture}(10,8)(0,0)\put(5,3){\circle{7}}
\put(3,2){$\scriptscriptstyle g$}\end{picture}$.

Define the genus $g(G)$ of a graph $G$ to be the sum of the genera $g(v)$
over the vertices $v$ of $G$ and the first Betti number of the graph
$G$. Then $g(G)$ equals the arithmetic genus of the curve $C$.

Elements of the link of a vertex $v$ are called the flags of the vertex,
and their number $n(v)$ is its valence: it equals the number of points in
the component $C(v)$ which map to a double point or marked point of
$C$. The flags of $G$ are the flags of its vertices, and the set of flags
is denoted $\Flag(G)$. The valence $n(G)$ of a graph $G$ is the number of
its legs.

The automorphism group $\Aut(G)$ of a dual graph is the subgroup of the
automorphism group of the underlying graph which fixes the legs, and
preserves the genera of the vertices.

\subsection{Stable curves}

We say that a pointed prestable curve is stable if it has no infinitesimal
automorphisms fixing the marked points, and that a dual graph is stable if
each vertex $v$ satisfies the condition $2(g(v)-1)+n(v)>0$; a pointed
prestable curve is stable if and only if the associated dual graph is
stable. The set of stable dual graphs $G$ of fixed genus $g(v)$ and valence
$n(v)$ is finite.

\section{Cycles in the moduli stack of prestable curves}

The moduli stack $\MMbar_{g,n}$ of prestable algebraic curves of arithmetic
genus $g$ with $n$ smooth marked points is a smooth algebraic stack of
dimension $3(g-1)+n$. (This is the stack denoted $\mathfrak{M}_{g,n}$ in
Behrend \cite{B}.) If $G$ is a dual graph of genus $g$ and valence $n$, let
$\MM(G)$ be the substack of $\MMbar_{g,n}$ consisting of the moduli of
pointed prestable curves with dual graph $G$. These substacks stratify
$\MMbar_{g,n}$, and the stratum $\MM(G)$ has codimension equal to the
number of edges of $G$.

We associate to a dual graph $G$ the cycle $\{G\}=[\MMbar(G)]$, which we will
frequently represent in formulas by the graph $G$ itself.

\subsection{Tautological classes and decorated dual graphs}
Denote by $\MMbar(G)$ the closure of $\MM(G)$ in $\MMbar_{g,n}$, by $i(G)$
the inclusion map $\MMbar(G)\hookrightarrow\MMbar_{g,n}$, and by $\tMM(G)$
the product
$$
\tMM(G) = \prod_{v\in\VERT(G)} \MMbar_{g(v),n(v)} .
$$
(This is the stack denoted $\mathfrak{M}(G)$ in Behrend \cite{B}.) The
following lemma is clear.
\begin{lemma} \label{tMM}
There is a natural action of $\Aut(G)$ on $\tMM(G)$, and a natural morphism
$$
\pi(G) : \tMM(G)\to\MMbar(G) ,
$$
identifying $\tMM(G)/\Aut(G)$ with $\MMbar(G)$.
\qed\end{lemma}

Let $\CCbar_{g,n} \to \MMbar_{g,n}$ be the universal curve over
$\MMbar_{g,n}$, whose fibre at $(C,z_i)$ is $C$. The projection
$\CCbar_{g,n}\to\MMbar_{g,n}$ has $n$ canonical sections
$\sigma_i:\MMbar_{g,n}\to\CCbar_{g,n}$, corresponding to the $n$
marked points. Let $\om_{\CCbar_{g,n}/\MMbar_{g,n}}$ be the relative
dualizing sheaf, and consider the line bundles
$\Om_i=\sigma_i^*\om_{\CCbar_{g,n}/\MMbar_{g,n}}$, and the associated
divisors $\Psi_i = c_1(\Om_i)$. (The line bundle $\Om_i$ has fibre
$T^*_{z_i}C$ at the prestable pointed curve $(C,z_1,\dots,z_n)$.)

\begin{definition}
A decorated dual graph $(G,\phi)$ is a dual graph $G$ together with a
function $\phi$ from the set of flags $\Flag(G)$ of $G$ to the natural
numbers $\{0,1,\dots\}$.
\end{definition}

We represent a decorated dual graph graphically by drawing at each flag $i$
of the graph $G$ $\phi(i)$ arrow-heads pointing towards the vertex to which
the flag $i$ is attached.

Given a decorated dual graph $(G,\phi)$, let $\{G,\phi\}$ be the cycle
$$
\{G,\phi\} = \frac{1}{|\Aut(G)|} \, i(G)_* \pi(G)_* \! \prod_{i\in\Flag(G)}
\Psi_i^{\phi(i)} .
$$
The codimension of $\{G,\phi\}$ equals the number of edges of $G$ plus the
number of arrow-heads. A simple example is the divisor $\Psi_1$ in
$\Mbar_{2,1}$, associated to the decorated dual graph
$$
\begin{picture}(10,30)(60,745)
\put( 65,774){\circle{6}}
\put( 63.5,772.5){$\scriptscriptstyle2$}
\put( 65,746){\vector(0,1){25}}
\end{picture}
$$

For each $N\ge0$, let $\Pi^g_{n,N}:\MMbar_{g,n+N} \to \MMbar_{g,n}$ be
the projection which maps the pointed prestable curve
$(C,z_1,\dots,z_n,z_{n+1},\dots,z_{n+N})$ to $(C,z_1,\dots,z_n)$.
\begin{proposition} \label{Pi}
Let $(G,\phi)$ be a decorated dual graph of genus $g$ and valence $n$, and
for $v$ a vertex of $G$, let $G_v$ be the dual graph obtained from $G$ by
attaching an additional leg to $G$ at the vertex $v$, with label $n+1$, and
let $\phi_v$ be the decoration of $G_v$ assigning $0$ to the new flag of
$G_v$. Then
$$
\bigl(\Pi^g_{n,1}\bigr)^{-1} \{G,\phi\} = \sum_{v\in\VERT(G)}
\frac{|\Aut(G_v)|}{|\Aut(G)|} \{G_v,\phi_v\} .
$$
\end{proposition}
\begin{proof}
It is easy to see that $\bigl( \Pi^g_{n,N} \bigr)^* \Psi_i=\Psi_i$, $1\le
i\le n$; thus, it suffices to prove the proposition in the undecorated
case. But in that case, it is true almost by definition: the additional
marked point $z_{n+1}$ must lie on one of the smooth components of $C$, and
each one gives rise to a term in the expansion of
$\bigl(\Pi^g_{n,1}\bigr)^{-1}\{G\}$; the rational factors take account of
the differing normalization of the cycles $\{G_v\}$.
\end{proof}

\section{Stabilization}

The open substack $\Mbar_{g,n}\subset\MMbar_{g,n}$ of moduli of stable
curves is a Deligne-Mumford stack, empty unless $2(g-1)+n>0$. Denote
the restrictions of the line bundles $\Om_i$ to $\Mbar_{g,n}$ by
$\om_i$, and their Chern classes $c_1(\om_i)$ by $\psi_i$. (These
classes are denoted $K_i$ by Mumford.)

There is a morphism of stacks $\stab : \MMbar_{g,n} \to \Mbar_{g,n}$,
called stabilization, which associates to a pointed prestable curve
$(C,z_1,\dots,z_n)$ the pointed stable curve obtained by contracting chains
of rational curves in $C$, that is, rational components of valence less
than $3$; by Proposition 3 of Behrend \cite{B}, the morphism $\stab$ is
flat.

If $G$ is a stable graph of genus $g$ and valence $n$, let $[G]$ be the
cycle in $\Mbar_{g,n}$ defined by the closed substack $\Mbar(G)$. More
generally, if $(G,\phi)$ is a stable decorated dual graph (that is, a
decorated dual graph such that $G$ is stable), define $[G,\phi]$ by the
same procedure as was used to define $\{G,\phi\}$ for general decorated
graphs, substituting the embedding
$i(G):\Mbar(G)\hookrightarrow\Mbar_{g,n}$ for
$i(G):\MMbar(G)\hookrightarrow\MMbar_{g,n}$, the Chern class $\psi_i$ for
$\Psi_i$, and the \'etale morphism
$$
\pi(G) : \prod_{v\in\VERT(G)} \Mbar_{g(v),n(v)} \to \Mbar(G)
$$
for the \'etale morphism $\pi(G):\tMM(G)\to\MMbar(G)$ of Lemma \ref{tMM}.

If $G$ is an undecorated stable graph, we have $\stab^*[G]=\{G\}$. The
relationship between $\stab^*[G,\phi]$ and $\{G,\phi\}$ is less immediate
when $\phi$ is non-trivial; for this reason, we wish to distinguish the
decorated graph associated to $[G,\phi]$ from that associated to
$\{G,\phi\}$. This we do by drawing the arrow-heads in a different
style. For example, the divisor $\psi_1$ in $\Mbar_{2,1}$ is represented by
the decorated stable graph
$$
\begin{picture}(10,30)(60,745)
\put( 65,774){\circle{6}}
\put( 63.5,772.5){$\scriptscriptstyle2$}
\put( 65,746){\lamsvector(0,1){25}}
\end{picture}
$$

The following proposition permits the expansion of the cycle
$\stab^*[G,\phi]$ as a linear combination of cycles $\{G',\phi'\}$.
\begin{proposition} \label{stab}
On $\MMbar_{g,n}$, we have
$$
\Psi_i = \stab^*\!\psi_i + 
\begin{picture}(30,35)(120,760)
\put(135,771){\line( 0,-1){ 15}}
\put(135,755){\circle*{4}}
\put(135,755){\line( 0,-1){ 15}}
\put(133,730){$i$}
\put(135,775){\circle{8}}
\put(133,774){$\scriptscriptstyle g$}
\put(132,778){\line(-3, 4){ 10}}
\put(138,778){\line( 3, 4){ 10}}
\put(129,790){$\dots$}
\put(120,793){$1$}
\put(133,793){$\widehat\imath$}
\put(144,793){$n$}
\end{picture}
$$
\vskip25pt
\end{proposition}
\begin{proof}
Denote the dual graph on the right-hand side of the above formula by $G_i$,
and introduce analogous dual graphs $G_j$, $1\le j\le n$, in which the role
of $i$ and $j$ is exchanged. The complement of $\Mbar_{g,n}$ in
$\MMbar_{g,n}$ is a divisor with components $\MMbar(G_j)$, $1\le j\le
n$. Since the restrictions of $\Psi_i$ and $\stab^*\!\psi_i$ to
$\Mbar_{g,n}$ are equal, we see that a formula of the general form
$$
\Psi_i = \stab^*\!\psi_i + k_1\{G_1\} + \dots + k_n\{G_n\}
$$
holds, where $(k_1,\dots,k_n)$ are certain rational numbers. These numbers
may be characterized by the formula
$$
h^*(\Psi_j-\stab^*\psi_j)=k_j c_1\bigl( N_{\MM(G_j)}\MMbar_{g,n} \bigr) ,
$$
where $h$ is the locally closed embedding
$h:\MM(G_j)\hookrightarrow\MMbar_{g,n}$, and $N_{\MM(G_j)}\MMbar_{g,n}$ is
the normal bundle of the embedding. Our task is to show that
$k_j=\delta_{ij}$; this follows from the following lemma.
\end{proof}

\begin{lemma}
If $i=j$, then $h^*(\Om_j\o\stab^*\om_j^{\Dual}) \cong
N_{\MM(G_j)}\MMbar_{g,n}$, otherwise it is a trivial line bundle.
\end{lemma}
\begin{proof}
Let $G$ be a dual graph with two vertices, of genus $g_1$ and $g_2$ and
valence $n_1$ and $n_2$, and one edge. Thus $g(G)=g_1+g_2$ and
$n(G)=n_1+n_2$; assume that $n(G)>0$, so that $\Aut(G)$ is trivial. The
codimension-one stratum $\MM(G)$ decomposes as a product
$$
\MM(G) \cong \MM_{g_1,n_1+1} \times \MM_{g_2,n_2+1} ,
$$
and we have the following formula for its normal bundle:
$$
N_{\MM(G)}\MMbar_{g,n} \cong \Om_{n_1+1}^{\Dual} \boxtimes
\Om_{n_2+1}^{\Dual} .
$$

On $\MM_{0,2}\cong B\GG_m$, we have the formula $\Psi_1+\Psi_2=0$,
reflecting the fact that the representations of $\GG_m$ on $T_0\P^1$
and $T_\infty\P^1$ are dual to each other. It follows that
\begin{align*}
h^*(\Psi_j-\psi_j) &=
\begin{picture}(35,45)(115,755)
\put(135,771){\line( 0,-1){ 14}}
\put(135,755){\circle*{4}}
\put(135,740){\vector( 0, 1){ 13}}
\put(133,730){$j$}
\put(135,775){\circle{8}}
\put(133,774){$\scriptscriptstyle g$}
\put(132,778){\line(-3, 4){ 10}}
\put(138,778){\line( 3, 4){ 10}}
\put(129,788){$\dots$}
\put(120,793){$1$}
\put(133,793){$\widehat\jmath$}
\put(144,793){$n$}
\end{picture} -
\begin{picture}(35,45)(120,755)
\put(135,755){\vector( 0,1){ 16}}
\put(135,755){\circle*{4}}
\put(135,740){\line( 0, 1){ 13}}
\put(133,730){$j$}
\put(135,775){\circle{8}}
\put(133,774){$\scriptscriptstyle g$}
\put(132,778){\line(-3, 4){ 10}}
\put(138,778){\line( 3, 4){ 10}}
\put(129,788){$\dots$}
\put(120,793){$1$}
\put(133,793){$\widehat\jmath$}
\put(144,793){$n$}
\end{picture} = \hbox{}
- \begin{picture}(35,45)(120,755)
\put(135,771){\vector( 0,-1){ 14}}
\put(135,755){\circle*{4}}
\put(135,740){\line( 0, 1){ 13}}
\put(133,730){$j$}
\put(135,775){\circle{8}}
\put(133,774){$\scriptscriptstyle g$}
\put(132,778){\line(-3, 4){ 10}}
\put(138,778){\line( 3, 4){ 10}}
\put(129,788){$\dots$}
\put(120,793){$1$}
\put(133,793){$\widehat\jmath$}
\put(144,793){$n$}
\end{picture} -
\begin{picture}(40,45)(115,755)
\put(135,755){\vector( 0,1){ 16}}
\put(135,755){\circle*{4}}
\put(135,740){\line( 0, 1){ 13}}
\put(133,730){$j$}
\put(135,775){\circle{8}}
\put(133,774){$\scriptscriptstyle g$}
\put(132,778){\line(-3, 4){ 10}}
\put(138,778){\line( 3, 4){ 10}}
\put(129,788){$\dots$}
\put(120,793){$1$}
\put(133,793){$\widehat\jmath$}
\put(144,793){$n$}
\end{picture} \\[30pt]
&= c_1\bigl( N_{\MM(G_j)}\MMbar_{g,n} \bigr) .
\end{align*}
On the other hand, for $i\ne j$, it is clear that
$h^*\bigl(\Om_i\o\om_i^{\Dual}\bigr)$ is canonically trivializable, showing
that $h^*(\Psi_i-\psi_i)=0$.
\end{proof}

For each $N\ge0$, let $\pi^g_{n,N}:\Mbar_{g,n+N} \to \Mbar_{g,n}$ be the
projection which maps the pointed stable curve
$(C,z_1,\dots,z_n,z_{n+1},\dots,z_{n+N})$ to the stabilization of
$(C,z_1,\dots,z_n)$. The analogue of Proposition \ref{Pi} only holds in
this case in the undecorated case.
\begin{proposition} \label{pi}
Let $G$ be a stable graph of genus $g$ and valence $n$, and for $v$ a
vertex of $G$, let $G_v$ be the dual graph obtained from $G$ by attaching
an additional leg to $G$ at the vertex $v$, with label $n+1$. Then
$$
\bigl(\pi^g_{n,1}\bigr)^{-1} [G] = \sum_{v\in\VERT(G)}
\frac{|\Aut(G_v)|}{|\Aut(G)|} [G_v] .
\qed$$
\end{proposition}

On the other hand, to calculate the pullback of the classes $\psi_i$ by
$\pi^g_{n,N}$, we need the following formula, which is the special case for
zero-dimensional $V$ of the more general Proposition \ref{psi-pull-back-V},
proved later in this article:
\begin{equation} \label{psi-pull-back}
\psi_i = \bigl( \pi^g_{n,1} \bigr)^* \psi_i +
\begin{picture}(35,40)(120,760)
\put(135,771){\line( 0,-1){ 15}}
\put(135,755){\circle*{4}}
\put(134,754){\line(-3,-4){ 10}}
\put(136,754){\line( 3,-4){ 10}}
\put(129,788){$\dots$}
\put(121,791){$\scriptstyle 1$}
\put(134,791){$\scriptstyle\hat\imath$}
\put(144,791){$\scriptstyle n$}
\put(135,775){\circle{8}}
\put(133,774){$\scriptscriptstyle g$}
\put(132,778){\line(-3, 4){ 8}}
\put(138,778){\line( 3, 4){ 8}}
\put(121,732){$\scriptstyle i$}
\put(140,732){$\scriptstyle n+1$}
\end{picture}
\end{equation}
\vskip30pt

Combining all of these results, we obtain the following result, by
induction on $N$.
\begin{theorem} \label{Stab}
$$
\begin{picture}(230,65)(-90,730)
\put(-90,760){$\displaystyle
\bigl( \Pi^g_{n,N} \bigr)^* \Psi_i = \bigl( \pi^g_{n,N}
\bigr)^* \psi_i + \sum_{\substack{I\coprod J=\{1,\dots,n+N\} \\
i\in I\subset\{i,n+1,\dots,n+N\}}}$}
\put(135,771){\line( 0,-1){ 15}}
\put(135,755){\circle*{4}}
\put(134,754){\line(-3,-4){ 10}}
\put(136,754){\line( 3,-4){ 10}}
\put(129,742){$\dots$}
\put(132,730){$I$}
\put(135,775){\circle{8}}
\put(133,774){$\scriptscriptstyle g$}
\put(132,778){\line(-3, 4){ 8}}
\put(138,778){\line( 3, 4){ 8}}
\put(129,788){$\dots$}
\put(132,793){$J$}
\end{picture}$$
\end{theorem}

\section{Tautological classes in genus $0$ and $1$}

Eq.\ \ref{psi-pull-back} may be used to give explicit formulas for the
classes $\psi_i$ on $\Mbar_{0,n}$ and $\Mbar_{1,n}$, which we now
recall. The moduli space $\Mbar_{0,3}$ is zero-dimensional, so that the
classes $\psi_i$ also vanish in it. Applying \eqref{psi-pull-back}
iteratively, we obtain a formula for $\psi_i$ on $\Mbar_{0,n}$: if
$j,k\in\{1,\dots,\hat\imath,\dots,n\}$ are distinct from $i$, then
\begin{gather} \label{psi-zero}
\psi_i=\sum_{\substack{I\coprod J=\{1,\dots,n\} \\ i\in I ; j,k\in J}}
\begin{picture}(35,35)(120,765)
\put(135,775){\line( 0,-1){ 20}}
\put(135,755){\line(3,-4){ 10}}
\put(135,755){\line(-3,-4){ 10}}
\put(129,742){$\dots$}
\put(132,730){$I$}
\put(135,775){\line(-3, 4){ 10}}
\put(135,775){\line( 3, 4){ 10}}
\put(129,788){$\dots$}
\put(133,793){$J$}
\end{picture} \\ \notag
\end{gather}

Since the line bundles $\EE$ and $\om_1$ on $\Mbar_{1,1}$ are isomorphic,
we see that $\psi_1 = \lambda_1$. There is one dual graph of codimension
$1$ in $\Mbar_{1,1}$, namely
$$
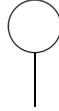
\begin{picture}(20,35)(70,750)
\put( 80,775){\circle{20}}
\put( 80,765){\line( 0,-1){20}}
\end{picture}
$$
with associated divisor $\delta_0$. The formula
$$
\lambda_1 = \frac{1}{12} \delta_0
$$
may be proved either by Grothendieck-Riemann-Roch or by consideration of
the explicit holomorphic section $\Delta$ (the discriminant) of
$\EE^{\o12}$, which has a simple pole at the divisor $\delta_0$. Combining
these two equations, we obtain the genus $1$ topological recursion relation
\begin{equation} \label{TRR-one}
\psi_1 = \frac{1}{12} \delta_0 .
\end{equation}
Applying \eqref{psi-pull-back} iteratively, we obtain a formula for
$\psi_i$ on $\Mbar_{1,n}$:
\begin{equation} \label{psi-one}
\psi_i = \frac{1}{12}
\begin{picture}(30,35)(120,760)
\put(135,755){\line(-3,-4){ 10}}
\put(135,755){\line( 3,-4){ 10}}
\put(129,742){$\dots$}
\put(122,732){$1$}
\put(142,732){$n$}
\put(135,765){\circle{20}}
\end{picture} +
\sum_{\substack{I\coprod J=\{1,\dots,n\} \\ i\in I}}
\begin{picture}(30,40)(120,760)
\put(135,773){\line( 0,-1){ 18}}
\put(135,755){\line(-3,-4){ 10}}
\put(135,755){\line( 3,-4){ 10}}
\put(129,742){$\dots$}
\put(132,730){$I$}
\put(135,775){\circle{4}}
\put(133,776){\line(-3, 4){ 9}}
\put(137,776){\line( 3, 4){ 9}}
\put(129,788){$\dots$}
\put(132,793){$J$}
\end{picture}
\end{equation}

\section{Gromov-Witten invariants}

\subsection{Prestable maps}

\begin{definition}
Let $V$ be a smooth projective variety. A pointed prestable map
$$\textstyle
(f:C\to V,z_1,\dots,z_n)
$$
is a pointed prestable curve together with an algebraic map $f:C\to
V$.
\end{definition}

Let $\NN_1(V)$ be the abelian group
$$
\NN_1(V) = \ZZ_1(V) / \textit{numerical equivalence} ,
$$
where $\ZZ_1(V)$ is the abelian group of $1$-cycles on $V$, and let
$\NE_1(V)$ be its sub-semigroup
$$
\NE_1(V) = \ZE_1(V) / \textit{numerical equivalence} ,
$$
where $\ZE_1(V)$ is the semigroup of effective $1$-cycles%
\footnote{Recall that two $1$-cycles $x$
and $y$ are numerically equivalent $x\equiv y$ when $x\*Z=y\*Z$ for any
Cartier divisor $Z$ on $V$.}.

\begin{lemma}[Proposition II.4.8, Koll\'ar \cite{Kollar}] \label{Mori}
If $V$ is a projective variety with K\"ah\-ler form $\om$, the set
$$
\{\beta\in\NE_1(V)\mid \om\cap\beta\le c\}
$$
is finite for each $c>0$.
\qed\end{lemma}

The dual graph of a prestable map with target $V$ is obtained from the dual
graph of the underlying pointed prestable curve by labelling each vertex
$v$ by the degree $\beta(v)\in\NE_1(V)$ of the restriction of $f$ to the
corresponding component $C(v)$ of $C$, that is, by the numerical
equivalence class of the $1$-cycle $f(C(v))$ in $\NE_1(V)$. The degree
$\beta(G)$ of the graph is the sum of these degrees over all vertices; it
equals the degree of the map $f$.

\subsection{Stable maps}

We say that a pointed prestable map is stable if it has no infinitesimal
automorphisms fixing the marked points, and that a dual graph is stable if
each vertex $v$ such that $\beta(v)=0$ satisfies the condition
$2(g(v)-1)+n(v)>0$; a pointed prestable map is stable if and only if the
associated dual graph is stable. Using Lemma \ref{Mori}, we see that the
set of stable dual graphs $G$ of fixed genus $g(v)$, valence $n(v)$ and
degree $\beta(v)$ is finite.

Behrend and Manin \cite{BM} show that Kontsevich's moduli stack
$\Mbar_{g,n}(V,\beta)$ of $n$-pointed stable maps of (arithmetic) genus $g$
and degree $\beta$ is a proper Deligne-Mumford stack (though not in general
smooth). It carries a forgetful map
$p:\Mbar_{g,n}(V,\beta)\to\MMbar_{g,n}$, obtained by discarding the map $f$
of a pointed stable map $(f:C\to V,z_1,\dots,z_n)$ --- the underlying
pointed curve $(C,z_1,\dots,z_n)$ of a stable map is only prestable, in
general.

The proof of the following result is similar to that of
Proposition~\ref{stab}.
\begin{proposition} \label{psi-pull-back-V}
Let $1\le i\le n$. On $\Mbar_{g,n+1}(V,\beta)$, we have the formula
$$
\begin{picture}(105,65)(45,730)
\put(45,760){$\pi^*p^*\Psi_i = p^*\Psi_i + \text{}$}
\put(135,771){\line( 0,-1){ 15}}
\put(135,755){\circle*{4}}
\put(141,753){$\scriptstyle 0$}
\put(134,754){\line(-3,-4){ 10}}
\put(136,754){\line( 3,-4){ 10}}
\put(129,788){$\dots$}
\put(121,791){$\scriptstyle 1$}
\put(134,791){$\scriptstyle\hat\imath$}
\put(144,791){$\scriptstyle n$}
\put(135,775){\circle{8}}
\put(133,774){$\scriptscriptstyle g$}
\put(142,772){$\scriptstyle \beta$}
\put(132,778){\line(-3, 4){ 8}}
\put(138,778){\line( 3, 4){ 8}}
\put(121,732){$\scriptstyle i$}
\put(140,732){$\scriptstyle n+1$}
\end{picture}
$$
\end{proposition}

\subsection{Virtual fundamental classes}
Let $\pi:\Mbar_{g,n+1}(V,\beta) \cong \Cbar_{g,n}(V,\beta) \to
\Mbar_{g,n}(V,\beta)$ be the universal curve, whose fibre over a pointed
stable map $(f:C\to V,z_1,\dots,z_n)$ is the curve $C$. Denote by
$f:\Cbar_{g,n}(V,\beta)\to V$ the universal stable map. If the sheaf
$R^1\pi_*f^*TV$ vanishes on $\Mbar_{g,n}(V,\beta)$, the Riemann-Roch
theorem predicts that the stack $\Mbar_{g,n}(V,\beta)$ is smooth, of
dimension
$$
\dim \Mbar_{g,n}(V,\beta) = (D-3)(1-g) + c_1(V)\cap\beta + n .
$$
This hypothesis is rarely true, and in any case only in genus
$0$. However, there is a cycle
$[\Mbar_{g,n}(V,\beta)]^\virt \in
A_{(D-3)(1-g)+c_1(V)\cap\beta+n}(\Mbar_{g,n}(V,\beta))$,
the virtual relative fundamental class, which stands in for
$[\Mbar_{g,n}(V,\beta)]$ in the obstructed case. Axiom I for Gromov-Witten
invariants (Behrend \cite{B}) says that, under the isomorphism
$\Mbar_{g,n}(V,0)\cong\Mbar_{g,n}\times V$, $[\Mbar_{g,n}(V,0)]^\virt$ is
identified with $[\Mbar_{g,n}]\times V$, while Axiom IV says that
$$
[\Mbar_{g,n+1}(V,\beta)]^\virt = \pi^*[\Mbar_{g,n}(V,\beta)]^\virt .
$$

\subsection{Gromov-Witten invariants}

In studying the Gromov-Witten invariants, it is convenient to work
with cohomology with coefficients in the Novikov ring $\Nov$ of $V$,
as adapted to the case of projective varieties (Lecture 4, Morrison
\cite{Morrison}; see also Section 2.1 of \cite{elliptic}). The Novikov
ring is
\begin{align*}
\Nov &= \Q[\NN_1(V)] \o_{\Q[\NE_1(V)]} \Q\[\NE_1(V)\] \\
&= \textstyle \bigl\{ a = \sum_{\beta\in\NN_1(V)} a_\beta q^\beta \mid
\text{ $\supp(a) \subset \beta_0+\NE_1(V)$ for some  $\beta_0\in\NN_1(V)$}
\bigr\} ,
\end{align*}
with product $q^{\beta_1}q^{\beta_2}=q^{\beta_1+\beta_2}$ and grading
$|q^\beta|=-2c_1(V)\cap\beta$. That the product is well-defined is shown by
Lemma \ref{Mori}.

We may now put these ingredients together to define the Gromov-Witten
invariants of $V$: if $x_1,\dots,x_n$ are cohomology classes of $V$,
\begin{multline*}
\< \tau_{k_1}(x_1) \dots \tau_{k_n}(x_n) \>_g \\ =
\sum_{\beta\in\NE_1(V)} q^\beta \int_{[\Mbar_{g,n}(V,\beta)]^\virt}
p^* \bigl( \Psi_1^{k_1} \dots \Psi_n^{k_n} \bigr) \cup \ev^*(x_1\boxtimes
\dots \boxtimes x_n) , 
\end{multline*}
where $\ev:\Mbar_{g,n}(V,\beta)\to V^n$ is evaluation at the marked points:
$$
\ev : { \textstyle (f:C\to V,z_1,\dots,z_n) } \mapsto \bigl(
f(z_1),\dots,f(z_n) \bigr) \in V^n .
$$
This is the $n$-point correlation function of two-dimensional topological
gravity with background the topological $\sigma$-model associated to
$V$. Observe that the grading of the Novikov ring is designed to ensure
that this correlator, considered as a map from $H^\bull(V^n,\Q)$ to $\Nov$,
is homogeneous, of degree
$$
2 \bigl\{ (3-D)(1-g) + \sum_{i=1}^n (k_i-1) \bigr\} .
$$

The following result is known as the string equation if $\om\in H^0(V,\Q)$,
and as the divisor equation if $\om\in H^2(V,\Q)$. It is an immediate
consequence of Axiom IV.
\begin{proposition} \label{string}
Let $\om\in H^i(V,\Q)$, $i\le 2$. Then
\begin{align*}
\<\om\tau_{k_1}(x_1)\dots\tau_{k_n}(x_n)\>_g &= {\textstyle \int_\beta\om }
\cdot \<\tau_{k_1}(x_1)\dots\tau_{k_n}(x_n)\>_g \\ &+ \sum_{i=1}^n
\<\tau_{k_1}(x_1)\dots\tau_{k_i-1}(\om\cup x_i)\dots\tau_{k_n}(x_n)\>_g ,
\end{align*}
except when $g=0$ and $n=2$: $\<\om x_1 x_2\>_0 = \int_V \om \cup x_1 \cup
x_2 + {\textstyle \int_\beta\om } \cdot \<x_1x_2\>_0$.  \qed
\end{proposition}

Another basic consequence of Axiom IV is the dilaton equation
$$
\<\tau_1(P)\tau_{k_1}(x_1)\dots\tau_{k_n}(x_n)\>_g = (2g-2+n)
\<\tau_{k_1}(x_1)\dots\tau_{k_n}(x_n)\>_g ,
$$
reflecting the fact that the relative dualizing sheaf
$\om_{\Mbar_{g,n+1}(V,\beta)/\Mbar_{g,n}(V,\beta)}$ has degree $2g-2$ along
the fibres, and that
$\om_{n+1}\o\om_{\Mbar_{g,n+1}(V,\beta)/\Mbar_{g,n}(V,\beta)}^{\Dual}$ is
the divisor of degree $n$ along the fibres associated to the marked
points. The string equation is equivalent to the equation $L_{-1}Z=0$ of
Eguchi et al.\ \cite{EHX}, where
$$
Z = \exp \biggl( \sum_{g=0}^\infty \lambda^{2g-2} \<\<~\>\>_g \biggr)
$$
while the divisor and dilaton equations together imply their equation
$L_0Z=0$.

\section{Topological recursion relations} \label{TRR}

We now combine the above results to show how to recover recursion relations
among Gromov-Witten invariants from relations in the Chow group of
$\Mbar_{g,n}$. To do this, we must turn our attention to Behrend's
construction of $[\Mbar_{g,n}(V,\beta)]^\virt$, since the axioms which he
states in \cite{B} are not quite adequate to the task.

If $G$ is a dual graph of genus $g=g(G)$ and valence $n=n(G)$, let
$\Mbar(G,V,\beta)$ be the fibred product
$\Mbar_{g,n}(V,\beta)\times_{\MMbar_{g,n}}\MMbar(G)$. If $G$ has edges
$\Edge(G)$, form the fibred product
$$\begin{CD}
\widetilde{\CM}(G,V,\beta) @>>>
\coprod\limits_{\beta=\sum_v\beta(v)}
\prod\limits_v \Mbar_{g(v),n(v)}(V,\beta(v)) \\
@VVV @VVV \\
V^{\Edge(G)} @>{\Delta}>> V^{\Edge(G)} \times V^{\Edge(G)}
\end{CD}$$
where $\Delta$ is the diagonal map.

There is an action of the group $\Aut(G)$ on $\widetilde{\CM}(G,V,\beta)$,
coming from compatible actions on $V^{\Edge(G)}$ and $\prod\limits_v
\Mbar_{g(v),n(v)}(V,\beta(v))$, and the quotient by this action is
naturally isomorphic to $\Mbar(G,V,\beta)$. Denote the quotient morphism
from $\widetilde{\CM}(G,V,\beta)$ to $\Mbar(G,V,\beta)$ by $\pi$.

The following proposition is the analogue of Lemma 10 of Behrend \cite{B}
for general dual graphs; it follows from Proposition 7.2 of \cite{BF} in
the same way as does Lemma 10. This result should probably be thought of as
an axiom for Gromov-Witten invariants, at the same level as Behrend's and
Manin's other axioms.
\begin{theorem} \label{Lemma10}
$$
i(G)^![\Mbar_{g,n}(V,\beta)]^\virt = \frac{1}{|\Aut(G)|} \, \pi_*\Delta^!
\!\!\!\!\!  \sum_{\beta=\sum_v\beta(v)} \bigotimes_v
[\Mbar_{g(v),n(v)}(V,\beta(v))]^\virt
$$
\end{theorem}

Taken together, Theorems \ref{Stab} and \ref{Lemma10} enable us to prove
the topological recursion relation \eqref{trr0}. Since $\Mbar_{0,3}$ is
$0$-dimensional, we have the equation $\psi_1=0$. It follows from
Proposition \ref{stab} that on $\MMbar_{0,3}$,
\begin{equation} \label{MM03}
\begin{picture}(120,40)(40,760)
\put( 90,760){$\Psi_1=$}
\put(135,773){\line( 0,-1){ 17}}
\put(135,755){\circle*{4}}
\put(135,755){\line( 0,-1){ 15}}
\put(132,730){$1$}
\put(135,775){\circle*{4}}
\put(133,776){\line(-3, 4){ 10}}
\put(137,776){\line( 3, 4){ 10}}
\put(120,793){$2$}
\put(144,793){$3$}
\end{picture}
\end{equation}
\vskip25pt
Call this dual graph $G$, and let $j$ be the inclusion of the stratum
$\MM(G)$ in $\MMbar_{0,3}$. We see that
\begin{multline*}
\< \tau_{k_1}(x_1)\tau_{k_2}(x_2)\tau_{k_3}(x_3) \>_0 \\
\begin{aligned}
{} &= \sum_\beta q^\beta \int_{[\Mbar_{0,3}(V,\beta)]^\virt}
\Psi_1^{k_1}\Psi^{k_2}\Psi^{k_3} \ev^*(x_1\boxtimes x_2\boxtimes x_3) \\
{} &= \sum_\beta q^\beta \int_{\Psi_1\cap[\Mbar_{0,3}(V,\beta)]^\virt}
\Psi_1^{k_1-1}\Psi^{k_2}\Psi^{k_3} \ev^*(x_1\boxtimes x_2\boxtimes x_3) \\
{} &= \sum_\beta q^\beta
\int_{j^![\Mbar_{0,3}(V,\beta)]^\virt} \Psi_1^{k_1-1}\Psi^{k_2}\Psi^{k_3}
\ev^*(x_1\boxtimes x_2\boxtimes x_3) \text{ by \eqref{MM03}} \\
{} &= \sum_a \< \tau_{k_1}(x_1) \gamma_a \>_0 \< \gamma^a
\tau_{k_2}(x_2) \tau_{k_3}(x_3) \>_0 \text{ by Theorem \ref{Lemma10}.}
\end{aligned}
\end{multline*}

This is not quite the same thing as \eqref{trr0}, but only a zeroth order
approximation to it, since it is an equation for
$\<\tau_{k_1}(x_1)\tau_{k_2}(x_2)\tau_{k_3}(x_3)\>_0$, and not for the
power series $\<\<\tau_{k_1}(x_1)\tau_{k_2}(x_2)\tau_{k_3}(x_3)\>\>_0$. To
prove \eqref{trr0}, we simply apply Theorem \ref{Stab} instead of
Proposition \ref{stab}.

We prove \eqref{trr1} in the same way. Since $\psi_1=\frac{1}{12}\delta_0$
on $\Mbar_{1,1}$, it follows from Proposition \ref{stab} that on
$\MMbar_{1,1}$,
$$
\Psi_1 = \frac{1}{12}
\begin{picture}(70,25)(65,760)
\put( 80,782){\circle{20}}
\put( 80,772){\circle*{4}}
\put( 80,772){\line( 0,-1){15}}
\put( 80,755){\circle*{4}}
\put( 80,753){\line( 0,-1){13}}
\put( 77,730){$1$}
\put(100,760){$+$}
\put(130,773){\line( 0,-1){ 17}}
\put(130,755){\circle*{4}}
\put(130,755){\line( 0,-1){ 15}}
\put(127,730){$1$}
\put(130,775){\circle{4}}
\end{picture}
$$
\vskip25pt\noindent
The proof of \eqref{trr1} now follows along the same lines as the
proof of \eqref{trr0}.

\section{Applications of the topological recursion relations}

The topological recursion relations \eqref{trr0} and \eqref{trr1} lead to
beautiful formulas of Dijkgraaf and Witten \cite{DW} which express
correlation functions on the large phase space in terms of correlation
functions on the small phase space. In this section, we follow their
proofs.
\begin{theorem}
Let $u^a$ be the power series $\<\<P\gamma^a\>\>_0$. Then
$$
\<\<\tau_{k_1}(x_1)\tau_{k_1}(x_2)\>\>_0 = \< \tau_{k_1}(x_1)
\tau_{k_2}(x_2) e^{u\*\gamma} \>_0 .
$$
\end{theorem}
\begin{proof}
Let $F_{k_1,k_2}(x_1,x_2)$ and $\Phi_{k_1,k_2}(x_1,x_2)$ be the left and
right-hand sides of the above equation. Denote equality modulo
$\{t^a_k\mid k>0\}$ by $a\sim b$. The proof relies on the following
consequence of the string equation \eqref{string},
\begin{equation} \label{small}
u^a\sim t^a_0 ,
\end{equation}
which implies that
$$
F_{k_1,k_2}(x_1,x_2) \sim \Phi_{k_1,k_2}(x_1,x_2) .
$$

We now calculate the derivatives of $F_{k_1,k_2}$ and $\Phi_{k_1,k_2}$ with
respect to $t^b_n$, $n>0$. On the other hand,
\begin{align*}
\frac{\p \Phi_{k_1,k_2}(x_1,x_2)}{\p t^b_n} &= \sum_a \frac{\p u^a}{\p
t^b_n} \< \gamma_a \tau_{k_1}(x_1) \tau_{k_2}(x_2) e^{u\*\gamma} \>_0
\\ &= \sum_a \<\<\tau_n(\gamma_b)P\gamma^a\>\>_0 \< \gamma_a
\tau_{k_1}(x_1) \tau_{k_2}(x_2) e^{u\*\gamma} \>_0 \\ &= \sum_{a,c}
\<\<\tau_{n-1}(\gamma_b)\gamma^c\>\>_0 \<\<\gamma_cP\gamma^a\>\>_0 \<
\gamma_a \tau_{k_1}(x_1) \tau_{k_2}(x_2) e^{u\*\gamma} \>_0 \\ &=
\sum_{a,c} \<\<\tau_{n-1}(\gamma_b)\gamma^c\>\>_0 \frac{\p u^a}{\p
t^c_0} \< \gamma_a \tau_{k_1}(x_1) \tau_{k_2}(x_2) e^{u\*\gamma} \>_0
\\ &= \<\<\tau_{n-1}(\gamma_b)\gamma^c\>\>_0 \frac{\p
\Phi_{k_1,k_2}(x_1,x_2)}{\p t^c_0} .
\end{align*}
On the other hand, we have by \eqref{trr0} that
\begin{align*}
\frac{\p F_{k_1,k_2}(x_1,x_2)}{\p t^b_n} &=
\<\<\tau_n(\gamma_b)\tau_{k_1}(x_1)\tau_{k_2}(x_2)\>\>_0 \\ &= \sum_a
\<\<\tau_{n-1}(\gamma_b)\gamma^a\>\>_0
\<\<\gamma_a\tau_{k_1}(x_1)\tau_{k_2}(x_2)\>\>_0 \\
&= \<\<\tau_{n-1}(\gamma_b)\gamma^a\>\>_0 \frac{\p
F_{k_1,k_2}(x_1,x_2)}{\p t^a_0} .
\end{align*}
Induction in the order of vanishing of
$F_{k_1,k_2}(x_1,x_2)-\Phi_{k_1,k_2}(x_1,x_2)$ in the variables $\{t^a_k\mid
k>0\}$ shows that the two power series are equal.
\end{proof}

The corresponding theorem in genus $1$, also due to Dijkgraaf and Witten
\cite{DW}, is even simpler, since it does not involve derivatives of the
potential.
\begin{theorem}
$$
\<\<~\>\>_1 = \< e^{u\*\gamma} \>_1 + \frac{1}{24} \log \det \left(
\frac{\p u^a}{\p t^b_0} \right)
$$
\end{theorem}
\begin{proof}
For the proof, it is useful to introduce the notation $M^a_b=\p u^a/\p
t^b_0$. We calculate that
\begin{align*}
\frac{\p M_a^b}{\p t^c_n} &=
\<\<\tau_n(\gamma_c)P\gamma_a\gamma^b\>\>_0
= \frac{\p}{\p t^a_0} \<\<\tau_n(\gamma_c)P\gamma^b\>\>_0 \\
&= \frac{\p}{\p t^a_0} \sum_d \<\<\tau_{n-1}(\gamma_c)\gamma^d\>\>_0
\<\<P\gamma_d\gamma^b\>\>_0 \\
&= \sum_d \<\<\tau_{n-1}(\gamma_c)\gamma^d\gamma_a\>\>_0 M_d^b
+ \sum_d \<\<\tau_{n-1}(\gamma_c)\gamma^d\>\>_0
\frac{\p M_a^b}{\p t^d_0}
\end{align*}
Multiplying by the matrix $M^{-1}$ and forming the trace, we see that
$$
\Tr \left( \frac{\p M}{\p t^c_n} M^{-1} \right) = \sum_a
\<\<\tau_{n-1}(\gamma_c)\gamma^a\gamma_a\>\>_0 + \sum_d
\<\<\tau_{n-1}(\gamma_c)\gamma^d\>\>_0 \Tr \left( \frac{\p M} {\p
t^d_0} M^{-1} \right) .
$$

Denote the left-hand side of the formula to be proved by $G$, and the
right-hand side by $\Gamma$. Eq.\ \eqref{small} shows that
$M^a_b\sim\delta^a_b$, so that $\det(M^a_b) \sim 1$, and hence $G \sim
\Gamma$.  Let us calculate the derivatives of $G$ and $\Gamma$ with respect
to $t^b_n$, $n>0$. On the one hand,
\begin{align*}
\frac{\p \Gamma}{\p t^c_n} &= \sum_a \frac{\p u^a}{\p t^c_n} \< \gamma_a
e^{u\*\gamma} \>_1 + \frac{1}{24} \Tr \Bigl( \frac{\p M}{\p t^c_n}
M^{-1} \Bigr) \\ &= \sum_a \<\<\tau_n(\gamma_c)P\gamma^a\>\>_0 \<
\gamma_a e^{u\*\gamma} \>_1 + \frac{1}{24} \Tr \Bigl( \frac{\p M}{\p
t^c_n} M^{-1} \Bigr) \\ &= \sum_{a,d}
\<\<\tau_{n-1}(\gamma_c)\gamma^d\>\>_0 \left(
\<\<\gamma_dP\gamma^a\>\>_0 \< \gamma_a e^{u\*\gamma} \>_0 +
\frac{1}{24} \Tr \Bigl( \frac{\p M} {\p t^d_0} M^{-1} \Bigr) \right) \\
& \quad {} + \frac{1}{24} \sum_a
\<\<\tau_{n-1}(\gamma_c)\gamma^a\gamma_a\>\>_0 \\ &=
\sum_a \<\<\tau_{n-1}(\gamma_c)\gamma^a\>\>_0 \frac{\p \Gamma}{\p t^a_0} +
\frac{1}{24} \sum_a \<\<\tau_{n-1}(\gamma_c)\gamma^a\gamma_a\>\>_0 .
\end{align*}
On the other hand, we have by \eqref{trr1} that
\begin{align*}
\frac{\p G}{\p t^c_n} &= \<\<\tau_n(\gamma_c)\>\>_1 = \sum_a
\<\<\tau_{n-1}(\gamma_c)\gamma^a\>\>_0 \<\<\gamma_a\>\>_1 +
\frac{1}{24} \sum_a \<\<\tau_{n-1}(\gamma_c)\gamma^a\gamma_a\>\>_0 \\
&= \sum_a \<\<\tau_{n-1}(\gamma_c)\gamma^a\>\>_0 \frac{\p G}{\p t^a_0} +
\frac{1}{24} \sum_a \<\<\tau_{n-1}(\gamma_c)\gamma^a\gamma_a\>\>_0 .
\end{align*}
Induction in the order of vanishing of $G-\Gamma$ in the variables
$\{t^a_k\mid k>0\}$ shows that the two power series are equal.
\end{proof}

\flushbottom

\part{The calculation of $\psi_1\psi_2$}

Mumford \cite{Mumford} has made a thorough analysis of the intersection
theory of the moduli space $\Mbar_2$. He shows that all of the cohomology
of $\Mbar_2$ is algebraic, and that it is generated by the boundary
divisors
$$
\delta_0 = \begin{picture}(40,25)(70,762)
\put( 80,782){\circle{5}}
\put( 80,780){\line( 0,-1){30}}
\put( 80,748){\circle{5}}
\end{picture}
\delta_1 = \begin{picture}(40,25)(60,762)
\put( 80,765){\circle{20}}
\put( 80,753){\circle{5}}
\end{picture} \\[10pt]
$$
There are two closed strata of codimension $2$ in $\Mbar_2$,
associated to the dual graphs
$$
\delta_{01} = \begin{picture}(80,15)(50,762)
\put( 62,765){\circle{5}}
\put( 64,765){\line(1,0){16}}
\put( 90,765){\circle{20}}
\end{picture}
\delta_{00} = \begin{picture}(60,15)(50,762)
\put( 70,765){\circle{20}}
\put( 90,765){\circle{20}}
\end{picture} \\[10pt]
$$
and in the cohomology ring of $\Mbar_2$, we have the formulas
\begin{equation} \label{M2}
\delta_0^2 = - 2 \delta_{01} + \frac{5}{3} \delta_{00} ,\quad
\delta_0\*\delta_1 = \delta_{01} ,\quad \delta_1^2 = - \frac{1}{12}
\delta_{01} .
\end{equation}

Faber~\cite{Faber} has extended Mumford's analysis to $\Mbar_{2,1}$. Once
more, all of the cohomology is algebraic, and is generated by the
tautological divisors:
$$
\begin{picture}(60,35)(35,740)
\put( 30,765){$\psi_1=$}
\put( 65,774){\circle{6}}
\put( 63.5,772.5){$\scriptscriptstyle2$}
\put( 65,746){\lamsvector(0,1){25}}
\end{picture}
\begin{picture}(60,45)(35,740)
\put( 30,765){$\delta_1=$}
\put( 65,783){\circle{5}}
\put( 65,766){\line(0,1){15}}
\put( 65,763){\circle{5}}
\put( 65,746){\line(0,1){15}}
\end{picture}
\begin{picture}(60,35)(100,722)
\put(100,747){$\delta_0=$}
\put(140,757){\circle{20}}
\put(140,745){\circle{4}}
\put(140,743){\line( 0,-1){15}}
\end{picture}
$$

\raggedbottom

We have the following formulas for intersection of pairs of divisors:
\begin{align*}
{} & \psi_1^2 = \frac{7}{5}
\begin{picture}(40,35)(60,775)
\put( 80,780){\line( 1, 1){ 13}}
\put( 80,780){\line(-1, 1){ 13}}
\put( 80,780){\line( 0,-1){ 18}}
\put( 66,796){\circle{5}}
\put( 94,796){\circle{5}}
\end{picture}
+ \frac{1}{5}
\begin{picture}(30,35)(65,757)
\put( 80,782){\circle{5}}
\put( 80,770){\circle{20}}
\put( 80,760){\line(0,-1){18}}
\end{picture}
- \frac{1}{120}
\begin{picture}(30,35)(65,757)
\put( 80,776){\circle{20}}
\put( 80,766){\line( 0,-1){9}}
\put( 80,755){\circle{5}}
\put( 80,753){\line( 0,-1){11}}
\end{picture}
+ \frac{13}{120}
\begin{picture}(30,35)(65,757)
\put( 80,770){\circle{20}}
\put( 80,760){\line( 3,-4){10}}
\put( 80,760){\line(-3,-4){12}}
\put( 91,745){\circle{5}}
\end{picture}
+ \frac{1}{120}
\begin{picture}(50,35)(55,760)
\put( 70,765){\circle{20}}
\put( 90,765){\circle{20}}
\put( 80,765){\line(0,-1){20}}
\end{picture} \\
{}& \psi_1\*\delta_1 = 2
\begin{picture}(35,25)(64,775)
\put( 80,780){\line( 1, 1){ 13}}
\put( 80,780){\line(-1, 1){ 13}}
\put( 80,780){\line( 0,-1){ 18}}
\put( 66,796){\circle{5}}
\put( 94,796){\circle{5}}
\end{picture}
+ \frac{1}{12}
\begin{picture}(47,35)(65,757)
\put( 80,770){\circle{20}}
\put( 80,760){\line( 3,-4){10}}
\put( 80,760){\line(-3,-4){12}}
\put( 91,745){\circle{5}}
\end{picture}
\psi_1\*\delta_0 =
2 \begin{picture}(29,35)(65,757)
\put( 80,782){\circle{5}}
\put( 80,770){\circle{20}}
\put( 80,760){\line(0,-1){17}}
\end{picture}
+ \begin{picture}(30,35)(65,757)
\put( 80,770){\circle{20}}
\put( 80,760){\line( 3,-4){10}}
\put( 80,760){\line(-3,-4){12}}
\put( 91,745){\circle{5}}
\end{picture}
+ \frac{1}{6}
\begin{picture}(50,35)(55,760)
\put( 70,765){\circle{20}}
\put( 90,765){\circle{20}}
\put( 80,765){\line(0,-1){20}}
\end{picture} \\
{}& \delta_1^2 =
- \frac{1}{12} \begin{picture}(30,35)(65,757)
\put( 80,776){\circle{20}}
\put( 80,766){\line( 0,-1){9}}
\put( 80,755){\circle{5}}
\put( 80,753){\line( 0,-1){11}}
\end{picture}
- \frac{1}{12}
\begin{picture}(60,35)(65,757)
\put( 80,770){\circle{20}}
\put( 80,760){\line( 3,-4){10}}
\put( 80,760){\line(-3,-4){12}}
\put( 91,745){\circle{5}}
\end{picture}
\delta_0\*\delta_1 =
\begin{picture}(30,35)(65,757)
\put( 80,776){\circle{20}}
\put( 80,766){\line( 0,-1){9}}
\put( 80,755){\circle{5}}
\put( 80,753){\line( 0,-1){11}}
\end{picture}
+
\begin{picture}(30,35)(65,757)
\put( 80,770){\circle{20}}
\put( 80,760){\line( 3,-4){10}}
\put( 80,760){\line(-3,-4){12}}
\put( 91,745){\circle{5}}
\end{picture} \\
{}& \delta_0^2 =
- 2
\begin{picture}(30,35)(65,757)
\put( 80,776){\circle{20}}
\put( 80,766){\line( 0,-1){9}}
\put( 80,755){\circle{5}}
\put( 80,753){\line( 0,-1){11}}
\end{picture}
- 2
\begin{picture}(30,35)(65,757)
\put( 80,770){\circle{20}}
\put( 80,760){\line( 3,-4){10}}
\put( 80,760){\line(-3,-4){12}}
\put( 91,745){\circle{5}}
\end{picture}
+ \frac{5}{3}
\begin{picture}(50,35)(55,760)
\put( 70,765){\circle{20}}
\put( 90,765){\circle{20}}
\put( 80,765){\line(0,-1){20}}
\end{picture} \\
\end{align*}
The first three of these formulas are proved by Faber~\cite{Faber}, using
the results of Section 8 of Mumford~\cite{Mumford}, while the last three
are obtained by pulling back the formulas \eqref{M2} by the morphism
$\pi^2_{1,1}:\Mbar_{2,2}\to\Mbar_{2,1}$.

The argument required to infer the first topological recursion relation
\eqref{trr2} from the above formula for $\psi_1^2$ are more elaborate than
the arguments in genus $0$ and $1$, since we must now deal with a quadratic
expression in divisors $\Psi_i$. Proposition \ref{stab} applied to the
class $\Psi_1$ on $\MMbar_{2,1}$, shows that
$$
\Psi_1^2 = \Psi_1 ( \stab^*\psi_1 ) +
\begin{picture}(20,30)(125,755)
\put(135,772){\line( 0,-1){ 16}}
\put(135,755){\circle*{4}}
\put(135,740){\vector( 0,1){ 13}}
\put(132,730){$1$}
\put(135,775){\circle{6}}
\put(133.5,773.5){$\scriptscriptstyle2$}
\end{picture}
= ( \stab^*\psi_1 )^2 +
\begin{picture}(20,30)(125,755)
\put(135,757){\lamsvector( 0,1){ 15}}
\put(135,755){\circle*{4}}
\put(135,740){\line( 0,1){ 13}}
\put(132,730){$1$}
\put(135,775){\circle{6}}
\put(133.5,773.5){$\scriptscriptstyle2$}
\end{picture} +
\begin{picture}(20,30)(150,755)
\put(160,772){\line( 0,-1){ 16}}
\put(160,755){\circle*{4}}
\put(160,740){\vector( 0,1){ 13}}
\put(157,730){$1$}
\put(160,775){\circle{6}}
\put(158.5,773.5){$\scriptscriptstyle2$}
\end{picture}
$$
\vskip27pt\noindent
A further application of Proposition \ref{stab} shows that
$$
\Psi_1^2 = ( \stab^*\psi_1 )^2 +
\begin{picture}(20,40)(125,755)
\put(135,772){\line( 0,-1){ 16}}
\put(135,755){\circle*{4}}
\put(135,740){\vector( 0,1){ 13}}
\put(132,730){$1$}
\put(135,775){\circle{6}}
\put(133.5,773.5){$\scriptscriptstyle2$}
\end{picture} +
\begin{picture}(20,40)(125,755)
\put(135,757){\vector( 0,1){ 15}}
\put(135,755){\circle*{4}}
\put(135,740){\line( 0,1){ 13}}
\put(132,730){$1$}
\put(135,775){\circle{6}}
\put(133.5,773.5){$\scriptscriptstyle2$}
\end{picture} -
\begin{picture}(20,40)(150,755)
\put(160,772){\line( 0,1){ 10}}
\put(160,770){\circle*{4}}
\put(160,757){\line( 0,1){ 11}}
\put(160,755){\circle*{4}}
\put(160,740){\line( 0,1){ 13}}
\put(157,730){$1$}
\put(160,785){\circle{6}}
\put(158.5,783.5){$\scriptscriptstyle2$}
\end{picture}
$$
\vskip27pt\noindent Combined with Mumford's formula \eqref{psi^2} for
$\psi_1^2$, this gives the topological recursion relation \eqref{trr2},
once we have replaced Proposition \ref{stab} by Theorem \ref{Stab}.

In order to prove the second topological recursion relation \eqref{trr22}
in genus $2$, we apply the expression \eqref{M22} for $\psi_1\psi_2$ in
terms of boundary classes, which we prove in the next few sections. The
topological recursion relation \eqref{trr22} follows by the familar pattern:
\begin{align*}
\Psi_1\Psi_2 = ( \stab^*\psi_1 ) \Psi_2 +
\begin{picture}(20,40)(125,760)
\put(135,785){\vector( 0,-1){ 12}}
\put(135,767){\line( 0,-1){ 13}}
\put(135,755){\circle*{4}}
\put(135,740){\line( 0,1){ 13}}
\put(132,730){$1$}
\put(132,788){$2$}
\put(135,770){\circle{6}}
\put(133.5,768.5){$\scriptscriptstyle2$}
\end{picture} &= \stab^*(\psi_1\psi_2) +
\begin{picture}(20,40)(125,760)
\put(135,785){\line( 0,-1){ 13}}
\put(135,768){\line( 0,-1){ 10}}
\put(135,770){\circle*{4}}
\put(135,740){\lamsvector( 0,1){ 12}}
\put(132,730){$1$}
\put(132,788){$2$}
\put(135,755){\circle{6}}
\put(133.5,753.5){$\scriptscriptstyle2$}
\end{picture} +
\begin{picture}(20,40)(125,760)
\put(135,785){\vector( 0,-1){ 12}}
\put(135,767){\line( 0,-1){ 13}}
\put(135,755){\circle*{4}}
\put(135,740){\line( 0,1){ 13}}
\put(132,730){$1$}
\put(132,788){$2$}
\put(135,770){\circle{6}}
\put(133.5,768.5){$\scriptscriptstyle2$}
\end{picture} \\[25pt]
&= \stab^*(\psi_1\psi_2) +
\begin{picture}(20,40)(125,760)
\put(135,785){\line( 0,-1){ 13}}
\put(135,768){\line( 0,-1){ 10}}
\put(135,770){\circle*{4}}
\put(135,740){\vector( 0,1){ 12}}
\put(132,730){$1$}
\put(132,788){$2$}
\put(135,755){\circle{6}}
\put(133.5,753.5){$\scriptscriptstyle2$}
\end{picture} -
\begin{picture}(20,40)(125,765)
\put(135,755){\circle*{4}}
\put(135,767){\circle{6}}
\put(133.5,765.5){$\scriptscriptstyle2$}
\put(135,779){\circle*{4}}
\put(135,764){\line(0,-1){10}}
\put(135,770){\line(0,1){10}}
\put(135,781){\line( 0,1){ 13}}
\put(135,753){\line( 0,-1){ 13}}
\put(132,730){$1$}
\put(132,796){$2$}
\end{picture} +
\begin{picture}(20,40)(125,760)
\put(135,785){\vector( 0,-1){ 12}}
\put(135,767){\line( 0,-1){ 13}}
\put(135,755){\circle*{4}}
\put(135,740){\line( 0,1){ 13}}
\put(132,730){$1$}
\put(132,788){$2$}
\put(135,770){\circle{6}}
\put(133.5,768.5){$\scriptscriptstyle2$}
\end{picture} \\[15pt]
\end{align*}

\flushbottom

\section{The Hodge polynomial of $\Mbar_{2,2}$}

In this section, we calculate the Hodge polynomials of $\Mbar_{2,2}$ and
$\Mbar_{2,2}/\SS_2$, where the symmetric group $\SS_2$ acts on on
$\Mbar_{2,2}$ through its action on the marked points of the stable curve
$(C,z_1,z_2)$.

Recall that the Serre characteristic of a quasi-projective variety, and
more generally, of a mixed Hodge module, is the Euler characteristic of its
cohomology in the Grothendieck group of mixed Hodge structures (see
\cite{II}). Since the Deligne-Mumford stacks $\Mbar_{2,2}$ and
$\Mbar_{2,2}/\SS_2$ are smooth and proper, their Serre characteristics
(that is, the Serre characteristic of the associated coarse moduli space)
determine their Hodge polynomials; thus, from now on, we will work only
with Serre characteristics.

An important role in our calculation will be played by the formula of
\cite{modular} for the $\SS_n$-equivariant Serre characteristics of the
compactifications $\Mbar_{g,n}$ in terms of the $\SS_n$-equivariant Serre
characteristics of the moduli stacks $\CM_{g,n}$:
\begin{multline} \label{Feynman}
\sum_{g=0}^\infty \hbar^{g-1} \sum_{2(g-1)+n>0}
\Serre^{\SS_n}(\Mbar_{g,n}) \\ = \sum_{k=0}^\infty \hbar^k D(h_k\circ
h_2) \* \sum_{g=0}^\infty \hbar^{g-1} \sum_{2(g-1)+n>0}
\Serre^{\SS_n}(\CM_{g,n}) .
\end{multline}

{}From \cite{II}, we know the equivariant Serre characteristics of the moduli
spaces $\CM_{0,n}$ and $\CM_{1,n}$. Thus, it only remains to calculate the
Serre characteristics of $\CM_{2,1}$, $\CM_{2,2}/\SS_2$ and $\CM_{2,2}$. It
will turn out that the calculation of the first two of these is relatively
straightforward, while to calculate the third, we must apply Faltings's
Eichler spectral sequence. Actually, for proving the topological recursion
relation in genus $2$, it is only necessary to calculate the Serre
characteristic of $\CM_{2,2}/\SS_2$, since the topological recursion
relation follows from a homological relation among $\SS_2$-invariant cycles
on $\Mbar_{2,2}$. We have chosen to include the calculation of
$\Serre(\CM_{2,2})$ since it implies that all of the cohomology of
$\Mbar_{2,2}$ is algebraic, a result which is of interest in itself.

\subsection{Mixed Hodge modules on $\CM_g$}

Let $\GSP(2g,\C)$ be the group of matrices
$$
\gamma = \begin{pmatrix} A&B\\C&D \end{pmatrix}
$$
such that $AB^T=BA^T$, $CD^T=DC^T$, and $AD^T-BC^T=\eta I_g$ for some
$\eta\in\C^\times$. The function $\eta:\GSP(2g,\C)\to\C^\times$
defines a character, the multiplier representation.

A maximal torus for $\GSP(2g,\C)$ is given by the diagonal matrices
$$
\gamma = \diag(a_1,\dots,a_g,d_1,\dots,d_g) ,
$$
such that $a_id_i=\eta$. The irreducible representations of $\GSP(2g,\C)$
have highest weight $$a_1^{\lambda_1}\dots
a_g^{\lambda_g}\eta^{n-|\lambda|},$$ where $\lambda_1\ge\dots\ge\lambda_g$
and $n\in\Z$; denote the associated representation by
$V_{\<\lambda\>}(n)$. For example, $V_{\<k\>}=V_{\<k\>}(0)$ is the $k$th
symmetric power of the contragradient representation of $\GSP(2g,\C)$ on
$\C^{2g}$, with weights $\{a_1/\eta,\dots,a_g/\eta,1/a_1,\dots,1/a_g\}$.

To an irreducible representation $V_{\<\lambda\>}(n)$ of $\GSP(2g,\C)$, we
associate a local system $\VV_{\<\lambda\>}(n)$ on $\CA_g$, underlying a
mixed Hodge module of weight $|\lambda|-2n$. We start with the fundamental
representation $V_{\<1\>}$, with associated mixed Hodge module
$\EE=R^1\bigl(\pi^g_{1,0}\bigr)_*\C$, where
$\pi^g_{1,0}:\Mbar_{g,1}\times_{\Mbar_g}\CM_g\to\CM_g$ is the universal
curve over $\CM_g$.  The representations $V_{\<\lambda\>}$ are the images
of certain idempotents in $(\C^{2g})^{\o|\lambda|}$; we apply the same
idempotents in $\EE^{\o|\lambda|}$ to obtain $\VV_{\<\lambda\>}$. Finally,
$\VV_{\<\lambda\>}(n)$ is obtained from $\VV_{\<\lambda\>}$ by tensoring
with $\C(n)$.

\subsection{The Serre characteristic of $\VV_{\<11\>}$}

Let $\Serre_\lambda$ be the Serre characteristic of the mixed Hodge module
$\VV_{\<\lambda\>}$ on $\CM_2$, and let $e_\lambda$ be its Euler
characteristic; note that $\Serre_0=q^3$, while $\Serre_\lambda=0$ if
$|\lambda|$ is odd. From the Leray spectral sequence and the techniques of
\cite{II}, we see that
\begin{align*}
\Serre(\CM_{2,1}) &= q^4 + q^3 , \\
\Serre^{\SS_2}(\CM_{2,2}) &= (q^5+q^4+\Serre_{11})s_2 + (q^4+\Serre_2)s_{11} , \\
\Serre^{\SS_3}(\CM_{2,3}) &= (q^6+q^5-q^4+q\Serre_{11}-\Serre_2)s_3
+ (q^5+q(\Serre_{11}+\Serre_2))s_{21} + (q+1)\Serre_2s_{111} .
\end{align*}
Applying \eqref{Feynman}, we obtain the Serre characteristics
of the compactified moduli spaces%
\footnote{These calculations were greatly aided by Stembridge's
\texttt{Maple} package \texttt{SF}.}:
\begin{align*}
\Serre(\Mbar_{2,1}) &= q^4+3q^3+5q^2+3q+1 , \\
\Serre^{\SS_2}(\Mbar_{2,2}) &= \bigl( q^5+5q^4+11q^3+11q^2+5q+2+\Serre_{11}
\bigr) s_2 \\ & \quad + \bigl( q^4+3q^3+3q^2+q+\Serre_2 \bigr) s_{11} , \\
\Serre^{\SS_3}(\Mbar_{2,3}) & =
(q^6+6q^5+20q^4+28q^3+20q^2+7q+2+(q+1)\Serre_{11}) s_3 \\ & \quad +
(3q^5+12q^4+19q^3+12q^2+4q+1+(q+1)(\Serre_{11}+\Serre_2))s_{21} \\ &
\quad + (q^3+(q+1)\Serre_2) s_{111} .
\end{align*}

\begin{proposition}
$\Serre_{11}=-1$
\end{proposition}
\begin{proof}
First, we show that the Euler characteristic $e_{11}=-1$. In the last
section of \cite{modular}, we prove that
$$
e(\CM_{0,6}/\SS_6) + e(\CM_{1,4}/\SS_4) + e(\CM_{2,2}/\SS_2) = 2 .
$$
But $e(\CM_{0,n}/\SS_n)=1$ for all $n\ge3$, $e(\CM_{1,4}/\SS_4)=0$ by
the results of \cite{II}, and $e(\CM_{2,2}/\SS_2)=2+e_{11}$.

The Poincar\'e polynomial of $\CM_3$ is $1+t^2+t^6$ (Looijenga
\cite{Looijenga}). The cohomology classes of degree $0$ and $2$ are
algebraic ($H^0(\CM_3,\C)$ is spanned by $1$ and $H^2(\CM_3,\C)$ is
spanned by $\lambda_1$), while the cohomology of degree $6$ must have
weight $w\in\{6,8,10,12\}$. It follows that
$\Serre(\CM_3)=q^6+q^5+q^{3-w/2}$. Applying \eqref{Feynman}, we see
that
$$
\Serre(\Mbar_3) = q^6 + 3q^5 + 7q^4 + 10q^3 + 7q^2 + 3q + 1 + q^{6-w/2} +
\Serre_{11} .
$$
Since $\Mbar_3$ satisfies Poincar\'e duality, we see that
$$
q^6 \Serre_{11}^\Dual = \Serre_{11} + q^{6-w/2} - q^{w/2} .
$$
On the other hand, since $\Mbar_{2,2}$ satisfies Poincar\'e duality,
we see that
$$
q^5 \Serre_{11}^\Dual = \Serre_{11} + 1 - q^5 .
$$
It follows that $e_{11}=11-w$, showing that $w=12$, and hence that
$\Serre_{11}=-1$.
\end{proof}

In the proof of this proposition, we could have used further results of
Looijenga, which show that $H^6(\CM_3,\Q)$ has weight $12$, but we thought
it interesting to present an alternate proof of this fact.

We conclude that
\begin{align} \label{Mbar22}
\Serre^{\SS_2}(\Mbar_{2,2}) &= \bigl( q^5+5q^4+11q^3+11q^2+5q+1\bigr) s_2
\\ & \quad + \bigl( q^4+3q^3+3q^2+q \bigr) s_{11} + \Serre_2s_{11} , \notag
\\ \Serre^{\SS_3}(\Mbar_{2,3}) & = (q^6+6q^5+20q^4+28q^3+20q^2+6q+1) s_3 \\
& \quad + (3q^5+12q^4+19q^3+12q^2+3q+)s_{21} + q^3 s_{111} \notag \\ &
\quad + \Serre_2(q+1)(s_{21}+s_{111}) . \notag
\end{align}
As we will see, $\Serre_2$ actually vanishes.

\subsection{Mixed Hodge modules on $\CA_g$}

Let $\CA_g$ be the moduli stack of completely polarized Abelian varieties
of genus $g$, and let $j_g:\CM_g\to\CA_g$ the Abel-Jacobi map. The mixed
Hodge module $\VV_{\<\lambda\>}(n)$ on $\CM_g$ is actually the pullback by
$j_g$ of a mixed Hodge modules of weight $2n-|\lambda|$ on $\CA_g$, which
we also denote by $\VV_{\<\lambda\>}(n)$, and which is defined in the same
way as the mixed Hodge complex $\VV_{\<\lambda\>}(n)$ on $\CM_g$, except
that we use the mixed Hodge module $\EE=R^1\pi_*\C$, where
$\pi:\CJ_g\to\CA_g$ is the universal Abelian variety over $\CA_g$, instead
of its pullback $j_g^*\EE\cong\bigl(\pi^g_{1,0})_*\C$.

In the special case of genus $2$, the Abel-Jacobi map
$j_2:\CM_2\hookrightarrow\CA_2$ is an open dense embedding, and the
associated Gysin long exact sequence for the mixed Hodge module
$\VV_{\<\lambda\>}$ reads
\begin{multline} \label{gysin}
\dots \to H^i_c(\CM_2,\VV_{\<\lambda\>}) \to
H_c^i(\CA_2,\VV_{\<\lambda\>}) \to
H_c^i(\CA_2\setminus\CM_2,\VV_{\<\lambda\>}) \\ \to
H_c^{i+1}(\CM_2,\VV_{\<\lambda\>}) \to H_c^{i+1}(\CA_2,\VV_{\<\lambda\>})
\to H_c^{i+1}(\CA_2\setminus\CM_2,\VV_{\<\lambda\>}) \to \dots
\end{multline}
Since $\VV_{\<\lambda\>}^\Dual\cong\VV_{\<\lambda\>}(|\lambda|)$,
Poincar\'e duality shows that
$$
H_c^i(\CA_2,\VV_{\<\lambda\>}) \cong
H^{6-i}(\CA_2,\VV_{\<\lambda\>})^{\Dual}(-|\lambda|) \cong
H^{6-i}(\SP(4,\Z),V_{\<\lambda\>})^{\Dual}(-|\lambda|)
$$
and
\begin{align*}
H_c^i(\CA_2\setminus\CM_2,\VV_{\<\lambda\>}) &\cong
H^{4-i}(\CA_2\setminus\CM_2,\VV_{\<\lambda\>})^{\Dual}(-|\lambda|) \\
&\cong H^{4-i}\bigl( \SS_2\wr\SL(2,\Z),
\Res^{\SP(4,\Z)}_{\SS_2\wr\SL(2,\Z)}V_{\<\lambda\>} \bigr)^{\Dual}(-|\lambda|) .
\end{align*}
Here, $\SS_2\wr\SL(2,\Z)$ is the subgroup of $\SP(4,\Z)$ which preserves
the standard decomposition of $\R^4$ into two symplectic vector spaces of
dimension $2$.

\raggedbottom

\subsection{The Eichler spectral sequence for $\GSP(4,\C)$}

Let $\Abar_g$ be a toroidal compactification of $\CA_g$, with
compactification divisor $D_\infty=\Abar_g\backslash\CA_g$ (see
Faltings-Chai \cite{FC}), and let $\pi:\Jbar_g\to\Abar_g$ be the
associated family of generalized Abelian varieties.

In the case that $g=2$, $\Abar_2$ is unique, and may be identified with
$\Mbar_2$ (Igusa \cite{Igusa}) in such a way that $\CA_2$ is identified
with the union of strata associated to the stable graphs
$\begin{picture}(10,10)(75,755) \put(80,760){\circle{6}}
\put(78.5,758.5){$\scriptscriptstyle2$}
\end{picture}$ and $\begin{picture}(45,10)(75,755) \put(
80,760){\circle{5}} \put( 82.5,760){\line(1,0){30}}
\put(115,760){\circle{5}}
\end{picture}$, and $D_\infty$ is identified with the stable graph
$\begin{picture}(25,15)(67,763)
\put( 80,765){\circle{20}}
\put( 80,753){\circle{5}}
\end{picture}$
\vskip12pt

There is an inclusion of $\GL(g,\C)$ in $\GSP(2g,\C)$, defined by
$$
A \mapsto \begin{pmatrix} A & 0 \\ 0 & A^{-1} \end{pmatrix} .
$$
Denote the irreducible representation of $\GL(g,\C)$ of highest weight
$a_1^{\lambda_1}\dots a_g^{\lambda_g}$ by $W_{(\lambda)}$. To such a
representation, let us associate a vector bundle $\CW_{(\lambda)}$ on
$\Abar_g$. Again, we start with the special case of the fundamental
representation $W_{(1)}\cong\C^g$: the associated vector bundle in this
case is $\CE=R^1\pi_*\CO$. The representations $W_{(\lambda)}$ are the
images of certain idempotents in $\bigl(\C^g\bigr){}^{\o|\lambda|}$
(different from those which occur in the symplectic case!); we may apply
the same idempotents in $\CE^{\o|\lambda|}$ to obtain $\CW_{(\lambda)}$.

In genus $2$, we will use the notation $\CE_k$ for the $k$-th symmetric
power $\CW_{(k)}$ of the Hodge bundle $\CE$, and $\lambda$ for its
determinant line bundle. Observe that
$\CW_{(k\ell)}\cong\lambda^\ell\o\CE_{k-\ell}$, and that
$\CW_{(k\ell)}^\Dual\cong\lambda^{-k}\o\CE_{k-\ell}$. We have the
identification $\CE_2\cong\Om^1_{\Abar_2}(\log D_\infty)$; taking the
determinant, we see that $\lambda^3\cong\om_{\Abar_2}(D_\infty)$, so that
Serre duality gives
\begin{align*}
& H^i_c(\Abar_2,\lambda^a\o\CE_b)^\Dual \cong
H^{6-i}(\Abar_2,\lambda^{3-a-b}\o\CE_b(D_\infty)) \text{ and } \\
& H^i_c(\Abar_2,\lambda^a\o\CE_b(-D_\infty))^\Dual \cong
H^{6-i}(\Abar_2,\lambda^{3-a-b}\o\CE_b) \text{ and }
\end{align*}

In Chapter VI of \cite{FC}, Faltings constructs a spectral sequence
converging to $H_c^\bull(\CA_g,\VV_{\<\lambda\>})$, which collapses at the
$E_1$-term, and whose $E_1$-term is a sum of cohomology groups of vector
bundles $\CW_{(\mu)}(-D_\infty)$. In genus $2$, his results amount to the
following theorem.
\begin{theorem}
The Hodge filtration on the cohomology groups
$H_c^\bull(\CA_2,\VV_{\<k\ell\>})$ has the form $0 \subset F^{k+\ell+3}
\subset F^{k+2} \subset F^{\ell+1} \subset F^0 =
H_c^\bull(\CA_2,\VV_{\<k\ell\>})$, where
\begin{align*}
F^0/F^{\ell+1} &\cong
H^\bull\bigl( \Abar_2 , \lambda^{-k}\o\CE_{k-\ell}(-D_\infty) \bigr) , \\
F^{\ell+1}/F^{k+2} &\cong
H^{\bull-1}\bigl( \Abar_2 , \lambda^{-k}\o\CE_{k+\ell+2}(-D_\infty) \bigr) , \\
F^{k+2}/F^{k+\ell+3} &\cong
H^{\bull-2}\bigl( \Abar_2 , \lambda^{1-\ell}\o\CE_{k+\ell+2}(-D_\infty) \bigr) , \\
F^{k+\ell+3} &\cong
H^{\bull-3}\bigl( \Abar_2 , \lambda^{\ell+3}\o\CE_{k-\ell}(-D_\infty) \bigr) .
\end{align*}
\end{theorem}
\begin{proof}
We show how this theorem follows from Theorem VI.5.5 of \cite{FC}, on
taking $g=2$ and highest weight $\lambda=a_1^ka_2^\ell$; this highest
weight arises because the dual of the associated mixed Hodge module is
$\VV_{\<k\ell\>}$.

The Weyl group of $\SP(4,\C)$ is generated by $\sigma_1$ and $\sigma_2$, where
\begin{align*}
& \sigma_1(a_1)=a_2, \sigma_1(a_2)=a_1, \sigma_1(d_1)=d_2 \text{ and }
\sigma_1(d_2)=d_1, \\
& \sigma_2(a_1)=a_1, \sigma_2(a_2)=d_2, \sigma_2(d_1)=d_1 \text{ and }
\sigma_2(d_2)=a_2.
\end{align*}
The character $\rho$ equals $a_1^2a_2\eta^{-3/2}$, the subset
$W^{\mathbf{M}}$ is
$\{e,\sigma_2,\sigma_2\sigma_1,\sigma_2\sigma_1\sigma_2\}$, and we have
\begin{align*}
e(a_1^ka_2^\ell\rho)/\rho &= a_1^ka_2^\ell, \\
\sigma_2(a_1^ka_2^\ell\rho)/\rho &= a_1^ka_2^{-\ell-2}\eta^{\ell+1} ,
\\ \sigma_2\sigma_1(a_1^ka_2^\ell\rho)/\rho &=
a_1^{\ell-1}a_2^{-k-3}\eta^{k+2} , \\
\sigma_2\sigma_1\sigma_2(a_1^ka_2^\ell\rho)/\rho &=
a_1^{-\ell-3}a_2^{-k-3}\eta^{k+\ell+3} .
\end{align*}
Taking the dual of the bundle with the highest weight associated to the
element of $W^{\mathbf{M}}$ with length $p$ and twisting by the divisor
$-D_\infty$, we obtain $E_1^{p\bull}$, and the theorem follows.
\end{proof}

\subsection{The mixed Hodge complex $\VV_{\<2\>}$}

\flushbottom

The Leray spectral sequence for $H^\bull_c(\CM_{2,2},\C)$ associated to the
projection $\pi^2_{0,2}:\CM_{2,2}\to\CM_2$ has $E_2$-term
$$
H_c^p\bigl(\CM_2,R^q\bigl(\pi^2_{0,2}\bigr){}_!\C \bigr) \cong
\begin{cases} H_c^p(\CM_2,\C)\o\C(-2) , & q=4 , \\
H_c^p(\CM_2,\EE\oplus\EE)\o\C(-1) = 0 , & q=3 , \\
H_c^p(\CM_2,\C)\o\C(-1) \oplus H^p(\CM_2,\EE^{\o2}) , & q=2 , \\
0 , & q<2 . \end{cases}
$$
Note that $\EE^{\o2} \cong \C(-1) \oplus \VV_{\<11\>} \oplus
\VV_{\<2\>}$. We now assemble some results on the cohomology of local
systems on $\CM_2$ which will help us to analyse this spectral sequence.

\begin{proposition}
The rational cohomological dimension $\cd_\Q\CM_2$ of $\CM_2$ is $3$.
\end{proposition}
\begin{proof}
The existence of an \'etale map $\CM_{0,6}\to\CM_2$, with Galois group
$\SS_6$ shows that $\cd_\Q \CM_2 = \cd_\Q \CM_{0,6}$.

If $n\ge3$, there is a fibration $\CM_{0,n+1}\to\CM_n$, with fibre
$\P^1\setminus\{z_1,\dots,z_n\}$; it follows from the Leray spectral
sequence that
$$\cd_\Q\CM_{0,n+1}=\cd_\Q\CM_{0,n}+\cd_\Q(\P^1\setminus\{z_1,\dots,z_n\}).$$
But $\cd_\Q(\P^1\setminus\{z_1,\dots,z_n\})=1$; since $\CM_{0,3}$ is
zero-dimensional, we see that $\cd_\Q\CM_{0,n}=n-3$.
\end{proof}

The rigidity theorem of Raghunathan \cite{Raghunathan} combined with the
Gysin exact sequence \eqref{gysin} together yield the following result
(which was pointed out to us by R. Hain).
\begin{proposition}
If $|\lambda|>0$, $H^0(\CM_2,\VV_{\<\lambda\>})$ and
$H^1(\CM_2,\VV_{\<\lambda\>})$ vanish.
\end{proposition}

Combining these observations, we see that the Leray spectral sequence for
$\pi^2_{2,0}:\CM_{2,2}\to\CM_2$ collapses, giving
$$
H_c^i(\CM_{2,2},\C) \cong \begin{cases} \C(-5) , & i=10 , \\
\C(-4)\oplus\C(-4) , & i=8 , \\
H_c^4(\CM_2,\VV_{\<11\>}) \oplus H_c^4(\CM_2,\VV_{\<2\>}) , & i=6, \\
H_c^3(\CM_2,\VV_{\<11\>}) \oplus H_c^3(\CM_2,\VV_{\<2\>}) , & i=5, \\ 0 , &
\text{otherwise.}
\qed\end{cases}
$$

\begin{lemma} \label{Serre_2}
$$
\Serre_2 = \bigl[ H^4_c(\CA_2,\VV_{\<2\>}) \bigr] - \bigl[
H^3_c(\CA_2,\VV_{\<2\>}) \bigr] + q
$$
\end{lemma}
\begin{proof}
By the Gysin exact sequence, we have
$$
\Serre_2 = \sum_i (-1)^i \bigl[ H^i_c(\CA_2,\VV_{\<2\>}) \bigr]
- \sum_i (-1)^i \bigl[
H^i_c(\CM_{1,1}\times\CM_{1,1},\VV_{\<2\>})^{\SS_2} \bigr] .
$$
By the branching rules for $\SL(2,\C)\times\SL(2,\C)\subset\SP(4,\C)$, the
restriction of $\VV_{\<2\>}$ to $\CM_{1,1}\times\CM_{1,1}$ is isomorphic to
$\VV_2\o\VV_0\oplus\VV_0\o\VV_2\oplus\VV_1\o\VV_1$, where $\VV_k$ is the
$k$th symmetric power of the Hodge local system on $\CM_{1,1}$. We have
$$
H^3_c(\CM_{1,1}\times\CM_{1,1},\VV_2\o\VV_0) \cong
H^3_c(\CM_{1,1}\times\CM_{1,1},\VV_0\o\VV_2) \cong \C(-1) ,
$$
and all other cohomology groups vanish. The exchange map on
$\CM_{1,1}\times\CM_{1,1}$ interchanges these two cohomology groups, and we
see that
$$
H^3_c(\CM_{1,1}\times\CM_{1,1},\VV_2\o\VV_0\oplus\VV_0\o\VV_2)^{\SS_2}
\cong \C(-1) ,
$$
proving the formula.
\end{proof}

We finally have enough information to deduce the following result.
\begin{proposition}
$\Serre_2=0$
\end{proposition}
\begin{proof}
We start by showing that $e_2=0$. We have seen that
$e(\CM_{2,3}/\SS_3)=1+e_{11}-e_2$. In the last section of \cite{modular},
we prove that
$$
e(\CM_{0,7}/\SS_7) + e(\CM_{1,5}/\SS_5) + e(\CM_{2,3}/\SS_3) +
e(\CM_{3,1}) = 6 .
$$
But $e(\CM_{0,n}/\SS_n)=1$ for all $n\ge3$, $e(\CM_{1,5}/\SS_5)=-1$ by
\cite{II}, and $e(\CM_{3,1})=6$ by Looijenga \cite{Looijenga}.

Since $H_c^3(\CM_2,\VV_{\<2\>})\subset H_c^5(\CM_{2,2},\C)$ and
$H_c^4(\CM_2,\VV_{\<2\>})\subset H_c^6(\CM_{2,2},\C)$, we see that the
weights of the cohomology groups $H_c^3(\CM_2,\VV_{\<2\>})$ lie between $0$
and $5$, while those of $H_c^4(\CM_2,\VV_{\<2\>})$ lie between $2$ and $6$.

Furthermore, the Eichler spectral sequence with $(k,\ell)=(2,0)$,
shows that mixed Hodge structure on $H^3_c(\CA_2,\VV_{\<2\>})$ has
$F$-weights in $\{0,1,4,5\}$, and that on $H^4_c(\CA_2,\VV_{\<2\>})$
has $F$-weights in $\{1,4,5\}$, so that
\begin{multline*}
H_c^3(\CA_2,\VV_{\<2\>}) \cong H_c^3(-)^{(0,0)} \oplus
H_c^3(-)^{(1,0)} \oplus H_c^3(-)^{(0,1)} \oplus H_c^3(-)^{(1,1)} \\
\oplus H_c^3(-)^{(4,0)} \oplus H_c^3(-)^{(0,4)} \oplus
H_c^3(-)^{(5,0)} \oplus H_c^3(-)^{(4,1)} \oplus H_c^3(-)^{(1,4)}
\oplus H_c^3(-)^{(0,5)}
\end{multline*}
and
$$
H_c^4(\CA_2,\VV_{\<2\>}) \cong H_c^4(-)^{(1,1)} \oplus H_c^4(-)^{(4,1)}
\oplus H_c^4(-)^{(1,4)} \oplus H_c^4(-)^{(5,1)} \oplus H_c^4(-)^{(1,5)} .
$$
Poincar\'e duality for the compactification $\Mbar_{2,2}$ implies that
$\Serre_2$ satisfies the functional equation $\Serre_2^\Dual = q^{-5}
\Serre_2$, or, taking into account Lemma \ref{Serre_2}, that
\begin{align*}
& H_c^3(-)^{(0,0)} = 0 , \quad
H_c^3(-)^{(1,0)} = H_c^3(-)^{(0,1)} = 0 , \\
& H_c^3(-)^{(1,1)} \cong H_c^4(-)^{(1,1)} \oplus \C(-1) , \quad
H_c^3(-)^{(4,0)} = H_c^3(-)^{(0,4)} = 0 , \\
& H_c^3(-)^{(5,0)} = H_c^3(-)^{(0,5)} = H_c^4(-)^{(5,1)} = H_c^4(-)^{(1,5)}
= 0 .
\end{align*}
We conclude that
$$
\Serre_2 = \bigl[ H_c^4(-)^{(4,1)} \oplus H_c^4(-)^{(1,4)} \bigr] - \bigl[
H_c^3(-)^{(4,1)} \oplus H_c^3(-)^{(1,4)} \bigr] .
$$
Since $\Serre_2$ is concentrated in weight $5$, it must be effective by
\eqref{Mbar22}. Since its Euler characteristic vanishes, it must be zero.
\end{proof}

\section{The rational cohomology ring of $\Mbar_{2,2}$}

There is a spectral sequence associated to the stratification of
$\Mbar_{2,2}$ by $\{\CM(G)\}$:
$$
E_2^{pq} = \bigoplus_{|\Edge(G)|=-q} H_c^{p+q}(\CM(G),\Q) \Rightarrow
H^\bull(\Mbar_{2,2},\Q) .
$$
This spectral sequence carries a mixed Hodge structure, compatible with the
mixed Hodge structure on $H^\bull(\Mbar_{2,2},\Q)$, and since all rational
cohomology of degree $k$ on $\Mbar_{2,2}$ has weight $k$, only classes of
weight $k$ in the degree $k$ summand of $E_2$ can possibly survive to
$E_\infty$. In degree $2$, there are six such classes, associated to the
decorated stable graphs
$$
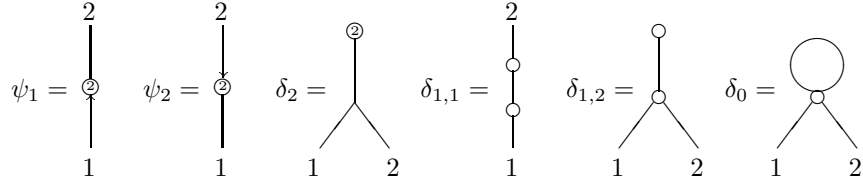
\begin{picture}(50,65)(35,740)
\put( 35,765){$\psi_1=$}
\put( 65,792){\line(0,-1){20}}
\put( 65,769){\circle{6}}
\put( 63.5,767.5){$\scriptscriptstyle2$}
\put( 65,746){\lamsvector(0,1){20}}
\put( 62,735){$1$}
\put( 62,794){$2$}
\end{picture}
\begin{picture}(50,65)(35,740)
\put( 35,765){$\psi_2=$}
\put( 65,792){\lamsvector(0,-1){19.5}}
\put( 65,769){\circle{6}}
\put( 63.5,767.5){$\scriptscriptstyle2$}
\put( 65,746){\line(0,1){20}}
\put( 62,735){$1$}
\put( 62,794){$2$}
\end{picture}
\begin{picture}(60,65)(105,720)
\put(105,745){$\delta_2=$}
\put(135,770){\circle{6}}
\put(133.5,768.5){$\scriptscriptstyle2$}
\put(135,767){\line( 0,-1){ 24}}
\put(135,743){\line(-3,-4){ 13}}
\put(135,743){\line( 3,-4){ 13}}
\put(117,715){$1$}
\put(147,715){$2$}
\end{picture}
\begin{picture}(50,65)(35,740)
\put( 30,765){$\delta_{1,1}=$}
\put( 65,792){\line(0,-1){12}}
\put( 65,777){\circle{5}}
\put( 65,775){\line(0,-1){12}}
\put( 65,760){\circle{5}}
\put( 65,746){\line(0,1){12}}
\put( 62,735){$1$}
\put( 62,794){$2$}
\end{picture}
\begin{picture}(60,65)(100,720)
\put(100,745){$\delta_{1,2}=$}
\put(135,770){\circle{5}}
\put(135,767){\line( 0,-1){ 20}}
\put(135,745){\circle{5}}
\put(133,743){\line(-3,-4){ 13}}
\put(137,743){\line( 3,-4){ 13}}
\put(117,715){$1$}
\put(147,715){$2$}
\end{picture}
\begin{picture}(50,45)(100,720)
\put(100,745){$\delta_0=$}
\put(135,757){\circle{20}}
\put(135,745){\circle{5}}
\put(133,743){\line(-3,-4){ 13}}
\put(137,743){\line( 3,-4){ 13}}
\put(117,715){$1$}
\put(147,715){$2$}
\end{picture}
$$
In degree $4$, there are fourteen such classes, associated to the decorated
stable graphs
$$
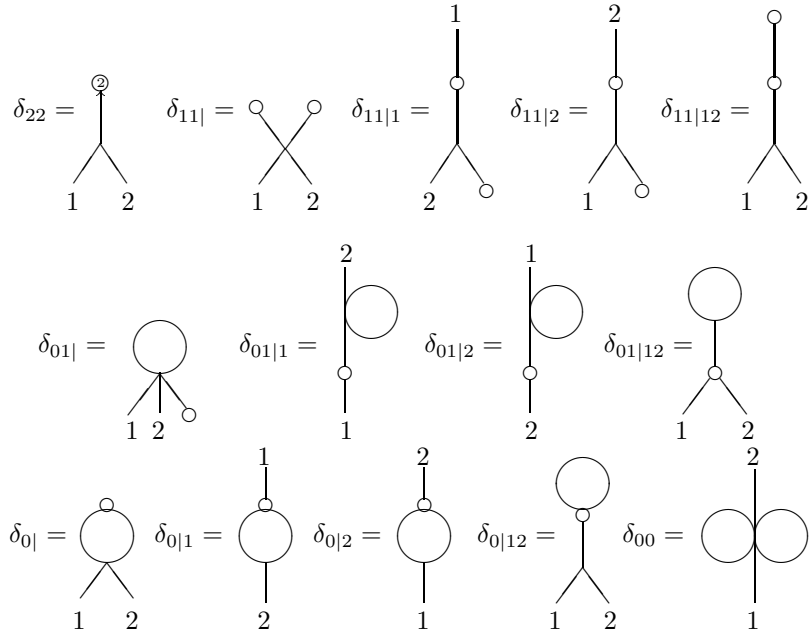
\begin{picture}(60,80)(45,722)
\put( 47,757){$\delta_{22}=$}
\put( 80,770){\circle{6}}
\put( 78.5,768.5){$\scriptscriptstyle2$}
\put( 80,747){\lamsvector( 0, 1){ 20}}
\put( 80,747){\line(-2,-3){ 10}}
\put( 80,747){\line( 2,-3){ 10}}
\put( 67,722){$1$}
\put( 88,722){$2$}
\end{picture}
\begin{picture}(65,80)(35,727)
\put( 35,762){$\delta_{11|}=$}
\put( 80,750){\line(-3,-4){ 10}}
\put( 80,750){\line( 3,-4){ 10}}
\put( 80,750){\line(-3, 4){ 10}}
\put( 80,750){\line( 3, 4){ 10}}
\put( 69,766){\circle{5}}
\put( 91,766){\circle{5}}
\put( 67,727){$1$}
\put( 88,727){$2$}
\end{picture}
\begin{picture}(60,80)(35,722)
\put( 40,757){$\delta_{11|1}=$}
\put( 77,793){$1$}
\put( 80,772){\line( 0, 1){ 18}}
\put( 80,770){\circle{5}}
\put( 80,767){\line( 0,-1){ 20}}
\put( 80,747){\line(-2,-3){ 10}}
\put( 80,747){\line( 2,-3){ 10}}
\put( 67,722){$2$}
\put( 91,729){\circle{5}}
\end{picture}
\begin{picture}(60,80)(35,722)
\put( 40,757){$\delta_{11|2}=$}
\put( 77,793){$2$}
\put( 80,772){\line( 0, 1){ 18}}
\put( 80,770){\circle{5}}
\put( 80,767){\line( 0,-1){ 20}}
\put( 80,747){\line(-2,-3){ 10}}
\put( 80,747){\line( 2,-3){ 10}}
\put( 67,722){$1$}
\put( 90,729){\circle{5}}
\end{picture}
\begin{picture}(60,80)(35,722)
\put( 37,757){$\delta_{11|12}=$}
\put( 80,795){\circle{5}}
\put( 80,772){\line( 0, 1){ 20}}
\put( 80,770){\circle{5}}
\put( 80,767){\line( 0,-1){ 20}}
\put( 80,747){\line(-2,-3){ 10}}
\put( 80,747){\line( 2,-3){ 10}}
\put( 67,722){$1$}
\put( 88,722){$2$}
\end{picture}
$$
$$
\begin{picture}(70,35)(35,745)
\put( 34,763){$\delta_{01|}=$}
\put( 80,765){\circle{20}}
\put( 80,755){\line(-3,-4){12}}
\put( 80,755){\line( 0,-1){15}}
\put( 80,755){\line( 3,-4){10}}
\put( 67,730){$1$}
\put( 77,730){$2$}
\put( 91,739){\circle{5}}
\end{picture}
\begin{picture}(70,65)( 90,735)
\put( 95,753){$\delta_{01|1}=$}
\put(135,785){\line( 0,-1){ 37}}
\put(145,768){\circle{20}}
\put(135,745){\circle{5}}
\put(135,742){\line( 0,-1){ 12}}
\put(133,720){$1$}
\put(133,787){$2$}
\end{picture}
\begin{picture}(70,65)( 90,735)
\put( 95,753){$\delta_{01|2}=$}
\put(135,785){\line( 0,-1){ 37}}
\put(145,768){\circle{20}}
\put(135,745){\circle{5}}
\put(135,742){\line( 0,-1){ 12}}
\put(133,720){$2$}
\put(133,787){$1$}
\end{picture}
\begin{picture}(70,65)( 90,735)
\put( 93,753){$\delta_{01|12}=$}
\put(135,775){\circle{20}}
\put(135,765){\line( 0,-1){ 17}}
\put(135,745){\circle{5}}
\put(133,743){\line(-3,-4){ 10}}
\put(137,743){\line( 3,-4){ 10}}
\put(120,720){$1$}
\put(145,720){$2$}
\end{picture} \\[15pt]
$$
$$
\begin{picture}(60,35)(40,762)
\put( 43,763){$\delta_{0|}=$}
\put( 80,765){\circle{20}}
\put( 80,777){\circle{5}}
\put( 80,755){\line(-3,-4){10}}
\put( 80,755){\line(3,-4){10}}
\put( 67,730){$1$}
\put( 87,730){$2$}
\end{picture}
\begin{picture}(60,35)(40,762)
\put( 38,763){$\delta_{0|1}=$}
\put( 80,765){\circle{20}}
\put( 80,777){\circle{5}}
\put( 80,755){\line(0,-1){15}}
\put( 80,779){\line(0,1){12}}
\put( 77,730){$2$}
\put( 77,792){$1$}
\end{picture}
\begin{picture}(60,35)(40,762)
\put( 38,763){$\delta_{0|2}=$}
\put( 80,765){\circle{20}}
\put( 80,777){\circle{5}}
\put( 80,755){\line(0,-1){15}}
\put( 80,779){\line(0,1){12}}
\put( 77,730){$1$}
\put( 77,792){$2$}
\end{picture}
\begin{picture}(60,35)(95,752)
\put( 95,753){$\delta_{0|12}=$}
\put(135,775){\circle{20}}
\put(135,763){\circle{5}}
\put(135,760){\line( 0,-1){ 17}}
\put(135,743){\line(-3,-4){ 10}}
\put(135,743){\line( 3,-4){ 10}}
\put(122,720){$1$}
\put(143,720){$2$}
\end{picture}
\begin{picture}(70,35)(25,762)
\put( 20,763){$\delta_{00}=$}
\put( 60,765){\circle{20}}
\put( 80,765){\circle{20}}
\put( 70,740){\line(0,1){50}}
\put( 67,730){$1$}
\put( 67,792){$2$}
\end{picture}
$$
\vskip35pt\noindent
As we saw in the last section, $h^2(\Mbar_{2,2})=6$ and
$h^4(\Mbar_{2,2})=14$; it follows that all of the classes which we have
enumerated must survive in $E_\infty$, with no relations among them. In
this way, we obtain the following result.
\begin{proposition}
Bases of the vector spaces $H^2(\Mbar_{2,2},\Q)$ and $H^4(\Mbar_{2,2},\Q)$
are given by the cycles associated to the above $6$ and $14$ decorated
stable graphs.
\end{proposition}

We now turn to the calculation of the intersection map
$$
H^2(\Mbar_{2,2},\Q) \o H^2(\Mbar_{2,2},\Q) \to H^4(\Mbar_{2,2},\Q)
$$
with respect to the given bases.
\begin{enumerate}
\item Six proper intersections among the divisors $\delta_0$,
$\delta_{1,1}$, $\delta_{1,2}$ and $\delta_2$:
$$\begin{tabular}{C|CCC}
& \delta_{1,1} & \delta_{1,2} & \delta_2 \\ \hline
\delta_0 & \delta_{01|1} + \delta_{01|2} & \delta_{01|} + \delta_{01|12}
& \delta_{0|12} \\
\delta_{1,1} & & \delta_{11|1} + \delta_{11|2} & 0 \\
\delta_{1,2} & & & \delta_{11|12}
\end{tabular}$$
\item The intersections $\psi_i\*\delta_{1,1}$ and $\psi_i\*\delta_{1,2}$
are calculated using \eqref{psi-one}:
\begin{align*}
\psi_1\*\delta_{1,1} &= \delta_{11|2} + \frac{1}{12} \delta_{01|2} , \\
\psi_2\*\delta_{1,1} &= \delta_{11|1} + \frac{1}{12} \delta_{01|1} , \\
\psi_1\*\delta_{1,2} &= 2 \delta_{11|} + \delta_{11|2} +
\delta_{11|12} + \frac{1}{12} \delta_{01|} , \\
\psi_2\*\delta_{1,2} &= 2 \delta_{11|} + \delta_{11|1} +
\delta_{11|12} + \frac{1}{12} \delta_{01|} .
\end{align*}
\item The intersection $\psi_1\*\delta_0$ is calculated by pulling back by
the morphism $\pi^2_{1,1}:\Mbar_{2,2}\to\Mbar_{2,1}$ the intersection
$\psi_1\*\delta_0$ in $\Mbar_{2,1}$. By \eqref{psi-pull-back},
$$
\bigl(\pi^2_{1,1}\bigr)^*\psi_1=\psi_1-\delta_2 ,
$$
showing that
$$
\bigl( \psi_1 - \delta_2 \bigr) \* \delta_0 = \psi_1\*\delta_0 -
\delta_{0|12} = \delta_{01|} + \delta_{01|2} + 2
\delta_{0|} + 2 \delta_{0|2} + \frac{1}{6} \delta_{00} .
$$
\item In the same way, the self-intersections $\psi_1^2$ and $\psi_2^2$ are
calculated by pullback from $\Mbar_{2,1}$:
\begin{align*}
\bigl(\pi^2_{1,1}\bigr)^* \bigl( \psi_1^2 \bigr) &= \bigl( \psi_1 -
\delta_2 \bigr)^2 = \psi_1^2 - \delta_{2,2} \\ &= \frac{1}{5} \bigl(
\delta_{0|} + \delta_{0|2} + 7 \delta_{11|} + 7 \delta_{11|2} \bigr) \\ &
\quad {} + \frac{1}{120} \bigl( \delta_{00} + 13 \delta_{01|} + 13
\delta_{01|2} - \delta_{01|1} - \delta_{01|12} \bigr) , \quad \text{and} \\
\psi_2^2 & = \frac{1}{5} \bigl( 7 \delta_{11|} + 7 \delta_{11|1}
+ \delta_{0|} + \delta_{0|1} \bigr) + \delta_{2,2} \\
& \quad {} + \frac{1}{120} \bigl( 13 \delta_{01|} + 13
\delta_{01|1} - \delta_{01|2} - \delta_{01|12} + \delta_{00} \bigr) .
\end{align*}
\item The self-intersection $\delta_0^2$ is calculated by pulling back the
formula \eqref{M2} for $\delta_0^2$ in $\Mbar_2$ by the morphism
$\pi^2_{0,2}:\Mbar_{2,2}\to\Mbar_2$:
$$
\delta_0^2 = - 2\delta_{01|} - 2\delta_{01|1} - 2\delta_{01|2} -
2\delta_{01|12} + \frac{5}{3} \delta_{00} .
$$
\item The self-intersections $\delta_{1,1}^2$, $\delta_{1,2}^2$ and
$\delta_2^2$ are calculated by the excess intersection formula, using the
explicit formulas for the normal bundles of these divisors and
\eqref{psi-one}:
\begin{gather*}
\delta_{1,1}\*\delta_{1,1} = -
\begin{picture}(16,40)(57,765)
\put( 65,792){\line(0,-1){12}}
\put( 65,777){\circle{5}}
\put( 65,774.5){\lamsvector(0,-1){12}}
\put( 65,760){\circle{5}}
\put( 65,746){\line(0,1){12}}
\put( 62,735){$1$}
\put( 62,794){$2$}
\end{picture} -
\begin{picture}(16,40)(57,765)
\put( 65,792){\line(0,-1){12}}
\put( 65,777){\circle{5}}
\put( 65,762.5){\lamsvector(0,1){12}}
\put( 65,760){\circle{5}}
\put( 65,746){\line(0,1){12}}
\put( 62,735){$1$}
\put( 62,794){$2$}
\end{picture}
\quad\quad
\delta_{1,2}\*\delta_{1,2} = -
\begin{picture}(50,40)(110,745)
\put(135,770){\circle{5}}
\put(135,767.5){\lamsvector( 0,-1){20}}
\put(135,745){\circle{5}}
\put(133,743){\line(-3,-4){ 13}}
\put(137,743){\line( 3,-4){ 13}}
\put(117,715){$1$}
\put(147,715){$2$}
\end{picture} -
\begin{picture}(50,40)(110,745)
\put(135,770){\circle{5}}
\put(135,747.5){\lamsvector( 0,1){20}}
\put(135,745){\circle{5}}
\put(133,743){\line(-3,-4){ 13}}
\put(137,743){\line( 3,-4){ 13}}
\put(117,715){$1$}
\put(147,715){$2$}
\end{picture} \\[20pt]
\delta_2^2 = -
\begin{picture}(30,40)(65,745)
\put( 80,770){\circle{6}}
\put( 78.5,768.5){$\scriptscriptstyle2$}
\put( 80,747){\lamsvector( 0, 1){ 20}}
\put( 80,747){\line(-2,-3){ 10}}
\put( 80,747){\line( 2,-3){ 10}}
\put( 67,722){$1$}
\put( 88,722){$2$}
\end{picture} -
\begin{picture}(30,40)(65,745)
\put( 80,770){\circle{6}}
\put( 78.5,768.5){$\scriptscriptstyle2$}
\put( 80,767){\lamsvector( 0, -1){ 20}}
\put( 80,747){\line(-2,-3){ 10}}
\put( 80,747){\line( 2,-3){ 10}}
\put( 67,722){$1$}
\put( 88,722){$2$}
\end{picture} = -
\begin{picture}(30,40)(65,745)
\put( 80,770){\circle{6}}
\put( 78.5,768.5){$\scriptscriptstyle2$}
\put( 80,747){\lamsvector( 0, 1){ 20}}
\put( 80,747){\line(-2,-3){ 10}}
\put( 80,747){\line( 2,-3){ 10}}
\put( 67,722){$1$}
\put( 88,722){$2$}
\end{picture} \\[5pt]
\end{gather*}
We find that $\delta_2^2 = - \delta_{22}$, and that
\begin{align*}
\delta_{1,1}^2 &= - \delta_{11|1} - \delta_{11|2} - \frac{1}{12}
(\delta_{01|1}+\delta_{01|2}) , \\
\delta_{1,2}^2 &= - \delta_{11|1} - \delta_{11|2} - \frac{1}{12}
(\delta_{01|}+\delta_{01|12}) .
\end{align*}
\end{enumerate}

We learn from these calculations that there are $7$ quadratic relations
among the six divisors $\psi_1$, $\psi_2$, $\delta_0$, $\delta_{1,1}$,
$\delta_{1,2}$, $\delta_2$. A basis for them is as follows:
\begin{align*}
& \delta_{1,1}\*(12\delta_{1,1}+12\delta_{1,2}+\delta_0) =
\delta_{1,2}\*(12\delta_{1,1}+12\delta_{1,2}+\delta_0) = 0 , \\
& \psi_1\*\delta_2 = \psi_2\*\delta_2 = 0 , \\
& \delta_{1,1}\*\delta_2 = \delta_{1,1}\*(\psi_1+\psi_2+\delta_{1,1}) = 0 , \\
& (\psi_1-\psi_2)\*(10\psi_1+10\psi_2-2\delta_{1,1}-12\delta_{1,2}-\delta_0)
= 0 .
\end{align*}

\section{The calculation of $\psi_1\*\psi_2$}

We now have all the data we need to prove the \eqref{M22}; in fact, by the
calculations of the last section, it will be a consequence of the following
result.
\begin{proposition} \label{basic}
$$
\psi_1\psi_2 = 3\,\psi_2^2 + \frac{1}{15} (3\,\psi_1-4\,\psi_2)
\delta_0 - \frac{1}{5} (8\,\psi_1+\psi_1) \delta_{1,1}
+ \frac{4}{5} (3\,\psi_1-4\,\psi_2) \delta_{1,2} - \frac{1}{3}
\delta_{0|1}
$$
\end{proposition}
\begin{proof}
Let $Q$ be the matrix of intersections between the $\binom{7}{2}=21$
quadratic monomials in the six divisors
$\{\psi_1,\psi_2,\delta_0,\delta_{1,1},\delta_{1,2},\delta_2\}$ and the
$\binom{8}{3}=56$ cubic monomials. An algorithm to calculate the entries of
this matrix has been implemented in \texttt{Maple} by Faber
\cite{Faber:Divisors} --- using this program and the code of Appendix A, we
may check that $Q$ has rank $14$, that its row space has basis
\begin{gather*}
\psi_i\psi_2 \quad \psi_i\delta_2 \quad \psi_i\delta_{1,1} \quad
\psi_i\delta_{1,2} \quad \psi_i\delta_0 \qquad (i=1,2) \\
\delta_2^2 \quad \delta_{1,2}\delta_2 \quad \delta_0\delta_2 \quad
\delta_0\delta_{1,1} \quad \delta_0\delta_{1,2} \quad \delta_0^2 ,
\end{gather*}
and that its column space has basis
\begin{gather*}
\psi_1\delta_{1,1}^2 \quad
\psi_1\delta_{1,1}\delta_{1,2} \quad
\psi_1\delta_{1,2}^2 \quad
\psi_1\delta_0^2 \\
\psi_1\psi_2\delta_0 \quad
\psi_1\psi_2\delta_{1,1} \quad
\psi_1\psi_2\delta_{1,2} \quad
\psi_1^2\delta_0 \quad
\psi_1\psi_2^2 \\
\delta_2\delta_0^2 \quad
\delta_2\delta_{1,2}^2 \quad
\delta_2^2\delta_{1,2} \quad
\psi_2\delta_{1,1}\delta_{1,2} \quad
\psi_2\delta_{1,2}^2 .
\end{gather*}
It follows that these sets are bases of $H^4(\Mbar_{2,2},\Q)$ and
$H^6(\Mbar_{2,2},\Q)$.

Of course, there are many such bases --- we have chosen these because it is
particularly easy to calculate the intersections of the cycle
$\delta_{0|1}$ with the above cubic monomials, since
$\delta_2\*\delta_{0|1}=0$ (these two cycles have empty intersection) and
$$
\psi_1\*\delta_{0|1} =
\begin{picture}(40,30)(60,762)
\put( 80,765){\circle{20}}
\put( 80,777){\circle{5}}
\put( 80,740){\lamsvector(0,1){15}}
\put( 80,779){\line(0,1){12}}
\put( 77,730){$2$}
\put( 77,792){$1$}
\end{picture} = 0 .
$$
\vskip30pt\noindent To calculate the two non-zero intersections
$\psi_1\delta_{1,2}^2\*\delta_{0|1}$ and
$\psi_1\delta_{1,1}\*\delta_{1,2}\*\delta_{0|1}$, we observe that
$$
\delta_{1,2} \* \delta_{0|1} =
\begin{picture}(25,35)(122,750)
\put(135,770){\circle{5}}
\put(135,767){\line( 0,-1){ 20}}
\put(135,745){\circle{5}}
\put(133,743){\line(-3,-4){ 13}}
\put(137,743){\line( 3,-4){ 13}}
\put(117,715){$1$}
\put(147,715){$2$}
\end{picture} \*
\begin{picture}(30,35)(65,765)
\put( 80,765){\circle{20}}
\put( 80,777){\circle{5}}
\put( 80,755){\line(0,-1){15}}
\put( 80,779){\line(0,1){12}}
\put( 77,730){$2$}
\put( 77,792){$1$}
\end{picture} =
\begin{picture}(40,35)(65,765)
\put( 80,765){\circle{20}}
\put( 80,755){\line(0,-1){15}}
\put( 80,775){\line(0,1){15}}
\put( 80,775){\line(1,0){18}}
\put(100,775){\circle{5}}
\put( 77,730){$2$}
\put( 77,792){$1$}
\end{picture} \text{ and }
\psi_1 \* \delta_{1,2} \* \delta_{0|1} =
\begin{picture}(40,35)(65,765)
\put( 80,765){\circle{20}}
\put( 80,755){\line(0,-1){15}}
\put( 80,790){\lamsvector(0,-1){15}}
\put( 80,775){\line(1,0){18}}
\put(100,775){\circle{5}}
\put( 77,730){$2$}
\put( 77,792){$1$}
\end{picture} =
\begin{picture}(40,40)(65,765)
\put( 80,765){\circle{20}}
\put( 80,755){\line(0,-1){15}}
\put( 80,805){\line(0,-1){30}}
\put( 80,790){\line(1,0){18}}
\put(100,790){\circle{5}}
\put( 77,730){$2$}
\put( 77,807){$1$}
\end{picture} \\[30pt]
$$
To calculate the intersections of this cycle with $\delta_{1,1}$ and
$\delta_{1,2}$, we use the excess intersection theorem, which shows
that
$$
\delta_{1,1} \*
\begin{picture}(40,35)(65,775)
\put( 80,765){\circle{20}}
\put( 80,755){\line(0,-1){15}}
\put( 80,805){\line(0,-1){30}}
\put( 80,790){\line(1,0){18}}
\put(100,790){\circle{5}}
\put( 77,730){$2$}
\put( 77,807){$1$}
\end{picture} = -
\begin{picture}(40,35)(65,775)
\put( 80,765){\circle{20}}
\put( 80,755){\line(0,-1){15}}
\put( 80,805){\lamsvector(0,-1){30}}
\put( 80,790){\line(1,0){18}}
\put(100,790){\circle{5}}
\put( 77,730){$2$}
\put( 77,807){$1$}
\end{picture} -
\begin{picture}(40,35)(65,775)
\put( 80,765){\circle{20}}
\put( 80,755){\line(0,-1){15}}
\put( 80,805){\line(0,-1){15}}
\put( 80,775){\lamsvector(0,1){15}}
\put( 80,790){\line(1,0){18}}
\put(100,790){\circle{5}}
\put( 77,730){$2$}
\put( 77,807){$1$}
\end{picture} = 0 \\[30pt]
$$
$$
\delta_{1,2} \*
\begin{picture}(40,50)(65,775)
\put( 80,765){\circle{20}}
\put( 80,755){\line(0,-1){15}}
\put( 80,805){\line(0,-1){30}}
\put( 80,790){\line(1,0){18}}
\put(100,790){\circle{5}}
\put( 77,730){$2$}
\put( 77,807){$1$}
\end{picture} = -
\begin{picture}(40,50)(65,775)
\put( 80,765){\circle{20}}
\put( 80,755){\line(0,-1){15}}
\put( 80,805){\line(0,-1){30}}
\put( 80,790){\lamsvector(1,0){17.5}}
\put(100,790){\circle{5}}
\put( 77,730){$2$}
\put( 77,807){$1$}
\end{picture} -
\begin{picture}(40,50)(65,775)
\put( 80,765){\circle{20}}
\put( 80,755){\line(0,-1){15}}
\put( 80,805){\line(0,-1){30}}
\put( 97.5,790){\lamsvector(-1,0){17.5}}
\put(100,790){\circle{5}}
\put( 77,730){$2$}
\put( 77,807){$1$}
\end{picture} = - \frac{1}{12}
\begin{picture}(60,50)(65,775)
\put( 80,765){\circle{20}}
\put( 80,755){\line(0,-1){15}}
\put( 80,805){\line(0,-1){30}}
\put( 80,790){\line(1,0){20}}
\put(110,790){\circle{20}}
\put( 77,730){$2$}
\put( 77,807){$1$}
\end{picture} = - \frac{1}{48} \\[45pt]
$$
The proposition now follows on solving the linear equation which expresses
the intersection vector of $\delta_{0|1}$ with our basis of cubic monomials
in terms of the intersection vectors for the basis of quadratic monomials.
\end{proof}

Observe that while our calculations do not respect the symmetry of
exchanging the labelling of the legs, the answer does --- as it must, since
$\psi_1\psi_2$ is invariant under this involution. This provides some
confirmation that we have performed the calculation correctly.

The proof of Proposition \ref{basic} also gives the following result.
\begin{theorem}
The six divisors $\psi_1$, $\psi_2$, $\delta_2$, $\delta_{1,1}$,
$\delta_{1,2}$ and $\delta_0$ generate the rational cohomology ring of
$\Mbar_{2,2}$.
\end{theorem}
\begin{proof}
The fact that the matrix $Q$ has rank $14$ implies that the subalgebra of
$H^\bull(\Mbar_{2,2},\Q)$ generated by these six divisors coincides with
$H^\bull(\Mbar_{2,2},\Q)$ up to degree $6$. The analogous intersection
matrix between linear and quartic monomials has rank $6$, completing the
proof.
\end{proof}

By contrast, $H^\bull(\Mbar_{2,3},\Q)$ is not generated by
$H^2(\Mbar_{2,3},\Q)$: the intersection matrix between quadratic and
quartic monomials in the $12$ divisors spanning $H^2(\Mbar_{2,3},\Q)$ has
rank $43$, while $h^4(\Mbar_{2,3})=44$. As we will show in a sequel to this
paper, the rational cohomology of $\Mbar_{2,3}$ is nevertheless all
algebraic.


\section*{Appendix A: \texttt{Maple} code used in proof of Lemma
\ref{basic}}

\begin{verbatim}
# Read in Faber's program for calculating intersection numbers.
read MgnF:
# generate the list of all [i_1,...,i_length]
# with 1<=i_1<=...<=i_L<=N
multi:=proc(N,L) local i,aux;
  aux:=proc(a) local j;
    seq([j,op(a)],j=1..op(1,a))
  end;
  if L=1 then [seq(i,i=1..N)] else map(aux,multi(N,L-1)) fi
end:
# Calculate the intersection matrix between monomials of degree
# d and 5-d in the divisors of Mbar_{2,2}. In this calculation,
# we omit the 3rd divisor (in Faber's ordering) kappa_1; in
# genus 2, this is a boundary divisor.
Intersection:=proc(d) local i,j,intersection,R,C;
  R:=subs(seq(i=i+1,i=n+1..6),multi(6,5-d));
  C:=subs(seq(i=i+1,i=n+1..6),multi(6,d));
  intersection:=array(1..nops(R),1..nops(C));
  for i from 1 to nops(R) do
    for j from 1 to nops(C) do
      intersection[i,j]:=mgn(2,[op(R[i]),op(C[j])])
    od od;
  RETURN(intersection)
end:
\end{verbatim}

\end{document}